\documentclass[11 pt, psamsfonts]{amsart}

\def\ignore#1{\relax}

\usepackage{heckearticle}

\usepackage{amssymb}

\input BoxedEPS

\SetRokickiEPSFSpecial
\HideDisplacementBoxes
\SetEPSFDirectory{./graphics/}

\ignore{
\SetTexturesEPSFSpecial
\HideDisplacementBoxes
\SetEPSFDirectory{:graphics:}
}

\begin{document}

\title[Hecke Algebras of type A]
{Iwahori-Hecke algebras of type A at roots of unity }

\author{Frederick M. Goodman and Hans Wenzl}

\address{Department of Mathemematics\\ University of Iowa\\ Iowa City, Iowa}

\email{goodman@math.uiowa.edu}

\address{Department of Mathematics\\ University of California\\ San Diego,
California}

\email{wenzl@brauer.ucsd.edu}

\thanks{We would like to thank the referee for an exceptionally careful
reading of a previous version of this paper, and for pointing the way towards
several improvements in the exposition.}

\maketitle
\bigskip

Let $H_n(q)$ be the Iwahori-Hecke algebra of type $A_{n-1}$ over a field $K$
of characteristic 0.
For each Young diagram $\la$, an $H_n(q)$ module
$S^\la$, called a  Specht module, was defined in [DJ].
The dimension of $S^\la$ does not depend on the choice of $q$,
and for $q$ not a root of unity, the $S^\la$'s provide a complete
set of irreducible $H_n(q)$-modules, up to isomorphism.
If $q$ is a primitive $l$-th root of unity, 
a complete set of simple $H_n(q)$-modules $D^\mu$
has been constructed in [DJ], where $\mu$ runs
through all Young diagrams with at most $l-1$ rows of
equal length. These modules are not well understood, but 
their dimensions can be computed 
 if one knows  the multiplicities $d_{\la\mu}$
of the  $D^\mu$ in a composition series of $S^\la$,
for all possible $\mu$'s and  $\la$'s.

A result known as Nakayama conjecture, or equivalently, in the
context of quantum groups, as linkage principle, gives some information
about the $d_{\la\mu}$'s. It says that $d_{\la\mu}\neq 0$ only
if $\mu+\rho$ is in the orbit of $\la +\rho$ with respect to 
an affine reflection group. Here $\rho = (k-1, k-2,\ ...,\ 1,0)$
for some $k\geq n$, and the reflection group is that
generated by the affine reflections in the hyperplanes $y_i-y_j=ml$,
$1 \le i < j \le k$ and 
$m\in\Z$. However, these orbits can become quite large and can
contain many diagrams $\mu$ for which $d_{\la\mu}=0$.

In this paper, we explore the use of {\em path idempotents} for the
Hecke algebra at roots of unity.  The path idempotents are defined via
the orthogonal representation, which is not well defined at roots of
unity.  Nevertheless, we show that certain sums of path idempotents are still 
well defined at roots of unity and can be used to derive information
about simple modules and decomposition numbers at roots of unity.

Our main technical result (Theorem \ref{Theorem 3.3.1}) concerns a (non-unital)
embedding of the ($k$-row quotient of) the Hecke algebra
$H_m(x^l)$ defined over the field of rational functions $K(x)$ into
the ($k$-row quotient of the) 
Hecke algebra $H_n(x)$, for certain $m$ and $n$. 
If $q$ is a primitive $l$-th root of unity, this gives 
an embedding of $KS_m$ into
$H_n(q)$. From this, we obtain the following applications:

(a) We show that $S^\la=D^\la$ for certain Young diagrams, namely for all 
diagrams $\la$ with $\leq k$ rows for which each component of $\la +\rho$
is divisible by $l$. This gives another proof  of special cases of results by
Dipper and James,  and of James and Mathas, 
for Specht modules.

(b) We obtain  estimates on the decomposition numbers
$d_{\la \mu}$  which provide rather good 
geometric information about the coefficients and in particular about
the set of $\la$ for which $d_{\la \mu} = 0$. These estimates
can be derived from the irreducibility criterion mentioned
above in fairly short order. (Such a derivation may be known  to experts,
but as far as we know it does not appear in the literature.)
The estimates certainly are easy consequences of
deep results by Soergel on Kazhdan-Lusztig polynomials and tilting
modules [S1,2].

Two algorithms for computing the $d_{\lambda \mu}$ have 
recently been discovered.  The first, proposed by Lascoux,
Leclerc, and Thibon [LLT] and proved by Ariki [Ar], computes
the decomposition numbers via the lower global crystal base
[K, Ha, MM] of the basic  $U_q(\hat{sl_l})$ module $L(\omega_0)$.
The second, due to Soergel [S1, S2] (see also [A2]) computes
decomposition numbers for tilting modules of $U_q(sl_k)$ via
certain Kazhdan-Lusztig polynomials for the affine Weyl group.
The proofs that the two methods give the decomposition numbers
[Ar, S2]  are quite deep, whereas the proofs of our estimates
are elementary.

It is also possible to give a more efficient version of the 
LLT algorithm which is a $q$-version of the algorithm for
our estimates.    This is  discussed in a separate paper [GW], where
we also show by a direct combinatorial method that the polynomial coeffcients 
of the global lower crystal base of $L(\omega_0)$ coincide with
the Kazhdan-Lusztig polynomials discussed by Soergel.

It would  be interesting
to get a better understanding between the connection of
our  principal technical result, the embedding of $H_m(x^l)$
into $H_n(x)$, and Lusztig's generalized Frobenius homomorphism 
for quantum groups. This
is discussed in greater
detail in Section  6.

\setcounter{section}{-1}
\section{Generic Algebras} 
This section contains some
preliminatries about finite dimensional algebras whose relations depend on a
parameter.  We suggest that the reader skip this section on the first reading,
noting however that it contains the definitions of some terms which
are used in the rest of the paper.  Recall that $K$ denotes a field of characteristic
0.

\subsection{Interpolation formulas}
Let $J$ be a finite set, and let $\{ a_j : j\in J\}$ be a set
of  variables.  We define for any
$l\in J$, and for any subset $S\subset J$ the Lagrange functions
$$\tilde P_{\ell}((a_j);y)=\prod_{j\neq \ell}\frac{(y-a_j)}
{(a_{l}-a_j)}\quad
{\rm and}\quad \tilde P_S=\sum_{\ell\in S} \tilde P_\ell;$$
these  are
rational functions in the variables $\{a_j\}$, and polynomials in the variable
$y$. Obviously, $\tilde P_{\ell}(a_i)=\delta_{i\ell}$.
One has the following well-known results regarding singularities of these
rational functions, and rational functional calculus of matrices.
\v\ni
\begin{lemma}\label{Lemma 0.1.1}
 The rational function
$\tilde P_S((a_j); y)$ has removable singularities at $a_i=a_j$ for
$i,j\in S$ as well as for $i,j \not\in S$ .
\end{lemma}
\v
\begin{proof} 
It follows from the definitions that we can write $\tilde P_S$ as a
rational function in the variables $\{a_j\}$ with denominator
$$\prod_{r, s \in S, r \ne s }(a_r-a_s)\ 
\prod_{s \in S, r \not\in S }(a_r-a_s).$$
The function $\tilde P_S$ is symmetric in the variables 
$\{a_j : j \in S\}$, while the denominator is antisymmetric in these
variables.  It follows that the numerator is also divisible by
$\prod_{r, s \in S, r \ne s }(a_r-a_s)$, so $\tilde P_S$  can also be written
as a rational function in the variables $\{a_j\}$ with denominator
$ \prod_{s \in S, r \not\in S }(a_r-a_s).\quad $
\end{proof}

\v\ni
\begin{corollary} \label{Corollary 0.1.2} Let $(S_i)_{i\in I}$ 
be a partition of the index set 
$J$.  Fix pairwise distinct elements $k_i\in K$,
for $i\in I$, and an index $\ell \in I$.  Then there exists a  polynomial
$P_{S_{\ell},(k_i)}(y)\in K[y]$ which is obtained from
$\tilde P_{S_{\ell}}((a_j);y)$ by substituting
$a_j=k_i$ when  $j\in S_i$. Moreover,
$P_{S_{\ell},(k_i)}(y)-1$ is divisible by $(y-k_{\ell})^s$, where
$s=|S_{\ell}|$.
\end{corollary}
\v
\ni
\begin{proof}  We obtain a well-defined polynomial
$P_{S_{\ell},(k_i)}$ by Lemma \ref{Lemma 0.1.1}. By construction, 
$\tilde P_{S_{\ell}}((a_j); a_r)=1$ for $a_r\in S_{\ell}$; hence
$\tilde P_{S_{\ell}}((a_j);y)-1$
is divisible by $\prod_{r\in S_{\ell}} (y-a_r)$.
The last statement follows from this. 
\end{proof}

\v\ni
\subsection{Eigenspaces} Fix once and for all an algebraic closure
$\overline{K(x)}$ of $K(x)$ and an algebraic closure $\overline K \subseteq
\overline{K(x)}$ of $K$.
Let $B$ be an $n\times n$ matrix with entries in $K(x)$
and with distinct eigenvalues $r_1, \dots, r_N$  in $\overline{K(x)}$,  with 
multiplicities $m_1, \dots, m_N$. Let $F = K(x,r_1, \dots, r_N)$. 
Recall that the generalized eigenspace of
$r_i$ is the kernel in $F^n$ of $(B-r_i1)^{m_i}$. Let $E(r_i)$ be the
projection onto this eigenspace, 
with kernel being the direct sum of the
other generalized eigenspaces; we call $E(r_i)$
the eigenprojection of $B$ for the eigenvalue $r_i$. 

We remark that it does not take much extra effort to deal with the case
that the matrix $B$ has its eigenvalues in an algebraic extension field
of $K(x)$, as we do here;  however, in our applications in the following sections, it
will never be necessary to consider such extension fields, so the reader can
safely ignore the extension fields here as well.

Let $\{v_j, 1 \le j \le n\}$ be a basis of $F^n$
consisting of generalized eigenvectors of $B$.
Partition  $J = \{i \in \N : 1 \le i \le n\}$
into subsets $S_i$, where $j\in S_i$ if $v_j$ is a generalized eigenvector
with eigenvalue $r_i$. Define the polynomials
$P_{S_{\ell}; (r_i)}$ for this partition, as in Corollary \ref{Corollary 0.1.2}.
\v\ni
\begin{lemma}
\label{Lemma 0.2.1} We have $E(r_{\ell})=P_{S_{\ell}; (r_i)}(B)$ for any eigenvalue
$r_\ell$ of $B$.
\end{lemma}
\v
\ni
\begin{proof}
By construction, the function $ P_{S_\ell;(r_i)}$ has a zero of
multiplicity $m_i$ at $r_i$ if $i\neq \ell$. By Corollary \ref{Corollary 0.1.2},
$ P_{S_\ell;(r_i)}-1$ has a zero of order $m_{\ell}$ at $r_{\ell}$.
Hence $ P_{S_\ell;(r_i)}(B)v_j=0$ for $v_j$ in the generalized
eigenspace of $r_i,\ i\neq \ell$, and 
$( P_{S_\ell;(r_i)}(B)-1)v_j=0$ for $v_j$ in the generalized
eigenspace of $r_{\ell}$.
\end{proof}
\v

\ni
\begin{definition}
Let $q\in K$. We call a matrix $B$ over $K(x)$ {\it evaluable} at $q$
(or just evaluable, if $q$ is understood)
if none of  its entries have poles  at $x=q$;  that is the entries of
the matrix lie in $K(x)_q$, the local ring of rational functions with
no poles at $q$.   The result of evaluating
the matrix at $q$  is denoted by $B(q)$. 
\end{definition}

In the following $B$ will denote an evaluable matrix. 
The eigenvalues $r_i$ of $B$ are algebraic integers over $K(x)_q$,
and 
the evaluation homomorphism from $K(x)_q$ to $K$ extends to a homomorphism
of $K(x)_q[r_1,\dots,r_N]$ to $\overline K$; denote the images of the $r_i$ by
$r_i(q)$. Let $F = K(x,r_1, \dots,r_N)$ and let 
$F_q$ denote the ring consisting of quotients of elements of
$K(x)_q[r_1,\dots,r_N]$  with denominators which are non-zero upon evaluation
at $q$; we refer to $F_q$ as the ring of evaluable elements of $F$. 
We call a matrix over $F$ evaluable if it has entries in $F_q$.

\v\ni
\begin{lemma}\label{Lemma 0.2.2} Suppose $B$ is an evaluable matrix over $K(x)$.
\begin{enumerate}
\item
The eigenvalues of $B(q)$ are given by $r_i(q)$.
\item
Suppose $B$ is invertible. Then 
$B\inv$ is evaluable if and only if $B(q)$ is invertible.
\end{enumerate}
\end{lemma}
\v\ni
\begin{proof} On the one hand, $\det(y \ident - B)$ maps to 
$\det(y \ident - B(q))$ under evaluation at $q$.  On the other hand,
$\det(y \ident - B) =  \prod (y-r_i)^{m_i}$ maps to  $\prod (y - r_i(q))^{m_i}$.
This shows part (a).  
If $B\inv$ is evaluable, then $B\inv(q)$  is an inverse for $B(q)$.
Conversely, if $B(q)$ is invertible, then $r_i(q) \ne 0$ for all $i$,
so $\det(B)(q) \ne 0$.  Then it follows from Cramer's rule that
$B\inv$ is evaluable. 
\end{proof}

\v\ni
\begin{proposition}
\label{Proposition 0.2.3}  Let $c\in K$ be an eigenvalue
of $B(q)$, and let $E(c)$ be its eigenprojection. Then 
the matrix $\displaystyle\sum_{\{\ell : r_\ell(q)=c\}} E(r_\ell)$ over $F$
is  evaluable  and its evaluation coincides with $E(c)$.
\end{proposition}
\v
\ni
\begin{proof}
For each eigenvalue $c$ of $B(q)$, let $T_c = \cup \{S_i :
r_i(q) = c\}$.  Then $\{T_c\}$ is a partition of $J$, and
$\sum_{\{\ell : r_\ell(q) = c\}} \tilde P_{S_\ell; (a_j)} = \tilde P_{T_c; (a_j)}$.
But $\tilde P_{T_c; (a_j)}$ can be written as
 a rational function in the variables $\{a_j\}$
with denominator  $\prod_{s \in T_c, r \not\in T_c}(a_r - a_s)$; therefore
$P_{T_c; (r_j)}$ lies in $F_q[y]$, and
$$\sum_{\{\ell : r_\ell(q)=c\}} E(r_\ell) =
\sum_{\{\ell : r_\ell(q)=c\}}  P_{S_\ell;(r_i)}(B) = P_{T_c; (r_j)}(B)
$$
is a matrix over $F_q$.  Furthermore,
$$ E(c)=  P_{T_c;(r_i(q))}(B(q))=
 P_{T_c;(r_i)}(B)(q)=
\left(\sum_{\{\ell : r_\ell(q)=c\}} E(r_\ell)\right)(q).$$
\end{proof}

\v\ni
\subsection{Evaluable elements of an algebra}
We apply the results of the previous section to 
the following set-up:
Let $A$ be a finite-dimensional $K(x)$-algebra with a  basis
$\{a_j : j\in J\}$, for which the structure coefficients are polynomials.

\begin{definition}
For  $q \in K$, 
we call an element $a=\sum_j s_j(x)a_j\in A$ {\em evaluable at $x=q$}
if none of the $s_j$  have a pole at $q$.  Denote by $A_q$ the set of 
evaluable elements of $A$.
\end{definition}
 
The set $A_q$ of evaluable elements coincides with the
span of $\{ a_j\}$ over $K (x)_q$, and is a $K (x)_q$-algebra.

Let $A(q)$ denote the $K$-algebra with a basis also denoted
by $\{a_j : j \in J\}$ whose structure coefficients are given by
evaluating the structure coefficients of $A$ at $q$.
Then $$A(q) \cong A_q \otimes_{K(x)_q} K,$$ where $K (x)_q$ acts on
$K$ by $f(x) \la = f(q) \la$.
The evaluation map 
$$a=\sum s_j(x)a_j \mapsto a(q)=\sum s_j(q)a_j,$$ or 
$a \mapsto a \otimes 1$,  
defines a $K$-algebra homomorphism from $A_q$ onto $A(q)$.

\begin{lemma} \label{Lemma 0.3.2} Let $A$ and $B$ be algebras over $K(x)$ with basis
$\{a_i\}, \{a_j\}$ having polynomial structure coefficients.  Let
$\varphi : A \rightarrow B$ be an algebra homomorphism such that
$\varphi(a_i) \in B_q$ for all $i$.
\begin{enumerate}
\item $\varphi(A_q) \subseteq B_q$.
\item  There is a unique $K$-algebra homomorphism $\tilde \varphi : A(q) \rightarrow
B(q)$ satisfying $\tilde \varphi (a(q)) = \varphi(a)(q)$ for $a \in A_q$.
\end{enumerate}
\end{lemma}

\begin{proof} Part (a) is evident.  The prescription $\tilde \varphi(a_i) =
\varphi(a_i)(q)$ determines the map \break $\tilde \varphi$ in part (b).
\end{proof}

\v\ni
\subsection{Evaluable representations}\ 
We will also apply the notion of evaluability  to
representations of $A$ with matrix coefficients in $K(x)$:

\begin{definition}
We say that a matrix representation $\Phi$ of $A$ is {\em evaluable} 
at $q$ if the matrices $\Phi(a_j)$ have coefficients
in $K(x)_q$ for all $j\in J$. An $A$-module $V$ is called {\em evaluable at $q$}
if it has a $K(x)$ basis with respect to which the basis elements
of $A$ act via $q$-evaluable matrices.
\end{definition}

{\em Observe that the class of the  $A(q)$
module may depend on the choice of basis.
 So we always assume for an evaluable $A$-module that a basis
has been fixed.}

\v\ni
\begin{lemma}
\label{Lemma 0.4.1} Let $\varphi$ be an evaluable representation of $A$ in
$\Mat_n(K(x))$.
\begin{enumerate}
\item
The restriction of $\varphi$ to the evaluable elements of $A$ is
a $K(x)_q$-algebra homomorphism into  $\Mat_n(K(x)_q)$.
\item
$\varphi$ induces a $K$-algebra homomorphism from 
$A(q)$ into $\Mat_n(K)$, determined by \break $\tilde \varphi(a(q)) = \tilde
\varphi(a)(q)$. 
\end{enumerate}
\end{lemma}
\v\ni
\begin{proof} This is a special case of Lemma \ref{Lemma 0.3.2}.\end{proof}
\v

Consider an evaluable representation $\varphi$ of $A$ in
$\Mat_n(K(x))$ and an evaluable element $a \in A$. 
Let $\{r_i\}$ be the spectrum in $\overline{K(x)}$
of left multiplication by $a$ in $A$, and put $F = K(x, r_1,\dots,r_N)$.  
The spectrum of $\varphi(a)$ is contained
in $\{r_i\}$.  Extend the evaluation homomorphisms to
$F_q$ as above.   Then the spectrum of $\varphi(a)(q)$ is contained in
the set $\{r_i(q)\}$.

\v\ni
\begin{lemma}
\label{Lemma 0.4.2} Let $\varphi$ be an evaluable representation of $A$.
\begin{enumerate}
\item
Let $a$ be an evaluable element of $A$.
Then for any eigenprojection $e'$ of $\varphi(a)(q)$
there exists an evaluable eigenprojection $e$ of $a$
such that its evaluation $\varphi(e)(q)$ coincides with $e'$.
\item
 Let $p'$ be an idempotent in $\varphi(A)(q)$. Then there exists
an evaluable idempotent $p\in A \otimes_{K(x)} F$, for
some finite algebraic extension $F$ of $K(x)$, 
such that $\varphi(p)(q)=p'$.
\item
Let $p'$ be an idempotent in $A(q)$.   Then there exists
an evaluable idempotent \break $p\in A \otimes_{K(x)} F$, for
some finite algebraic extension $F$ of $K(x)$, 
such that $p(q)=p'$.
\end{enumerate}
\end{lemma}
\v
\ni
\begin{proof} Let $c$ be the eigenvalue corresponding to $e'$. 
Then, by Proposition \ref{Proposition 0.2.3}, the eigenprojection $e$ onto the 
direct sum of generalized eigenspaces corresponding to the eigenvalues
$r_i$ with $r_i(q)=c$ acts via an evaluable matrix in the left
regular representation, i.e. it is an
evaluable element (in $A \otimes_{K(x)} F$). 
Hence also the matrix $\varphi(e)$ is evaluable and $\varphi(e)(q)$
 coincides with the eigenprojection $e'$ of $\varphi(a)(q)$,
again by Proposition \ref{Proposition 0.2.3}.

To show (b), let $a'\in A(q)$, say $a'=\sum_j s_jb_j$
with $s_j\in K$,  be such that $\varphi(a')=p'$.
Let $a=\sum_j s_jb_j\in A$, where $(b_j)$ is the corresponding
basis for $A$ and the $s_j$'s are constant rational functions.
Then $p'$ is an eigenprojection of $\varphi(a)(q)=p'$. Hence there
exists an evaluable idempotent $p\in A\otimes_{K(x)} F$ 
for an appropriate extension field $F$, such that $\varphi(p)(q)=p'$,
by part (a).

Point (c) is the special case of (b) where the representation $\varphi$ is
the left regular representation.  
\end{proof}
\v\ni
\subsection{Rank vectors} 
We shall henceforth assume that 
the algebra $A$ is a semisimple algebra over
$K(x)$ which is a direct sum of full matrix algebras over $K(x)$:
$$A=\bigoplus_{\la\in \Lambda} A_\la .$$
The examples we shall deal with, the
Hecke algebras of type $A$, satisfy this assumption.

\def\Tr{{\rm Tr}}

 Let $B$ be a semisimple algebra over a field $K$
which is isomorphic to the direct sum
$\bigoplus_\la B_\la$ of full matrix algebras $B_\la$.
Let $p$ be an idempotent in $B$.
The rank ${\bf r}(p)$ is the vector $(r(p)_\la)_\la$,
where $r(p)_\la = \Tr_\la(p)$ with $\Tr_\la$ the usual trace on
End($V_\la$) and $V_\la$ a simple $B_\la$ module.

In particular one can apply this notion of the rank vector to
idempotents in $A$.

The algebra $A(q)$ may no longer be semisimple. Let
$\overline{A(q)}$ be its maximum semisimple quotient,
with $\overline{A(q)}\cong \bigoplus_\mu \overline{A(q)}_\mu$ and 
$\overline{A(q)}_\mu$ a full matrix algebra over $K$ for each index $\mu$.
Let $D^\mu$ denote the simple $\overline{A(q)}_\mu$ module.
We define the (reduced)
rank  vector ${\bf r}_q(p')$ of an idempotent
$p'\in A(q)$ with respect to the maximum semisimple quotient $\overline{A(q)}$.

\v\ni 
\subsection{Evaluable modules} Recall that an $A$ module $W$ is called evaluable if
$W$ has a basis
with respect to which the matrices representing the basis elements
of $A$ are evaluable at $x=q$.  An evaluable $A$ module $W$ ``restricts" to
an $A(q)$ module $W_q$; as a set, $W_q$ is the $K$-linear span of 
the distinguished basis of $W$; thus $W = W_q \otimes_{K} K(x)$.

 Assume that $A$ has a complete set
of $q$-evaluable simple modules $(S^\la)$.

Consider a minimal idempotent
in $\overline{A(q)}$ which is in the simple component $\overline{A(q)}_\mu$ of
$\overline{A(q)}$. Such an idempotent  can be lifted to an idempotent $p'_\mu$
in $A(q)$, and two such liftings differ by a nilpotent element of the
radical of $A(q)$.

As $S^\la$ is evaluable, $p'_\mu$ also acts on $S^\la_q$.
Let $d_{\la\mu}=\Tr_{S^\la_q}(p_\mu')$ be its rank in this representation;
this is independent of the choice of the lifting since the difference 
of two liftings is nilpotent.

These numbers can be interpreted in terms of a composition series
of the $A(q)$ modules $S^\la_q$. More generally,
let $W$ be a $q$-evaluable $A$ module whose restriction to  an $A(q)$ module
has a composition series of the form
$$0=V_0\subset V_1\subset\ ...\ \subset V_k=W_q;\eqno(*)$$
here the $V_i$'s are $A(q)$ modules 
such that the factors $V_i/V_{i-1}$ are simple. 
\v\ni
\begin{proposition}
\label{Proposition 0.6.1} 
 Let $W$ be an evaluable module, and let
$p \in A$ be an evaluable idempotent.
\begin{enumerate}
\item
 $\Tr_W(p)=\Tr_{W_q}(p(q))$.
\item
 Suppose that $p(q)=p'_\mu$. Then ${\bf r}(p)=(d_{\la\mu})_\la$.
\item
If $p(q)=p'_\mu$, then $\Tr_W(p)$ is the number of factors in
the composition series $(*)$ which are isomorphic to $D^\mu$.
\item
 If $W\cong \oplus_\la m_\la S^\la$ as an $A$ module,
then exactly $\sum_\la m_\la d_{\la\mu}$ factors in $(*)$ are isomorphic
to $D^\mu$.
\end{enumerate}
\end{proposition}
\v
\ni
\begin{proof}
As $p$ is an idempotent, it acts via a matrix on $W$
whose diagonal entries (which are rational functions) add up
to a constant integer. Evaluating these functions at $q$ and summing
them up produces the same integer. This proves (a). Claim (b) follows
from (a) and the definition of $d_{\la\mu}$.

To prove (c), observe that $\Tr_{W_q}(p'_\mu)
=\sum_i \Tr_{V_{i+1}/V_i}(p'_\mu)$. 
As $V_{i+1}/V_i$ is simple and $p'_\mu$ is a lifting of 
a  minimal idempotent in
$\overline{A(q)}$, it follows that
$\Tr_{V_{i+1}/V_i}(p'_\mu)$ is either 0 or 1, depending
on whether $V_{i+1}/V_i$ is a simple $\overline{A(q)}_\mu$ module or not.
Hence $\Tr(p'_\mu)$ is equal to the number of factors isomorphic
to $D^\mu$. The statement now follows from (a).

Statement (d) follows from (b), (c) and $\Tr_W(p)
=\sum_\la m_\la \Tr_{S^\la}(p)$. 
\end{proof}
\v\ni
\begin{corollary}
\label{Corollary 0.6.2} If $p$ is an evaluable idempotent in $A$, we have
${\bf r}(p)=D\,{\bf r}_q(p(q))$, where $D=(d_{\la\mu})$.
\end{corollary}
\v\ni
\subsection{Existence of idempotents}
\begin{proposition}
\label{Proposition 0.7} Let $(\pi_\la)_\la$ be  not necessarily evaluable
representations of $A$ such that for each $\la$,
 $\pi_\la$ is equivalent over $K(x)$ to the 
representation on $S^\la$. 
Let $a$ be an evaluable element of $A$.  Then there exists an evaluable
idempotent $p$ in $A \otimes_{K(x)} F$, for some finite algebraic extension
$F$ of $K(x)$,  with the following property: 
whenever $\pi_\la(a)$ is an evaluable matrix, and its 
evaluation  $\pi_\la(a)(q)$ is a rank $d_\la$ idempotent, then 
$p$ acts as a rank $d_\la$ idempotent  on $S^\la$.
\end{proposition}
\v
\ni
\begin{proof} Consider $a$ acting via left multiplication on $A$.
Let $\{r_i\}$ be the spectrum of $a$ in $\overline{K(x)}$, 
 $F = K(x, r_1, \dots, r_N)$, and  let $F_q$ denote the ring of 
$q$-evaluable elements of $F$.
Let $p$ be the eigenprojection of $a$ projecting on the span
of generalized eigenspaces of $a$ belonging to the eigenvalues $r_i$
for which $r_i(q)=1$.   By Proposition \ref{Proposition 0.2.3}, $p$ is evaluable, and
$p=P(a)$, where $P \in F_q[y]$.

Suppose that $\pi_\la(a)(q)$ is a rank $d_\la$ idempotent.
Note that $\pi_\la(p) = \pi_\la(P(a)) = P(\pi_\la(a))$ is an
evaluable idempotent.  Furthermore its evaluation
$\pi_\la(p)(q) = \pi_\la(P(a))(q) = P(\pi_\la(a)(q))$ is the
eigenprojection of $\pi_\la(a(q))$ corresponding to its eigenvalue 1,
i.e. it is equal to   $\pi_\la(a)(q)$, by assumption.

Finally,
$$
\Tr_{S^\la}(p) = \Tr(\pi_\la(p)) = \Tr(\pi_\la(p)(q)) = \Tr(\pi_\la(a)(q)) = d_\la,
$$
where the first equation follows from the equivalence of the representations, 
the second equation results
from Proposition \ref{Proposition 0.6.1}(a) -- applied to the algebra generated by $p$,
the third  follows from the equality $\pi_\la(p)(q) = \pi_\la(a)(q)$
observed above, and the last from the assumption on $\pi_\la(a)(q)$. 
\end{proof}
\v
\ni
\begin{remark}
We mention again that in the applications of Lemma \ref{Lemma 0.4.2} and Proposition 
\ref{Proposition 0.7} in the following sections, we can work entirely over $K(x)$,
and need not pass to an extension field.  The idempotents obtained using these
results will lie in the algebra $A$ itself rather than in some $A \otimes_{K(x)} F$.
\end{remark}

\section{Hecke algebras}
\subsection{Definition of Hecke algebras} Let $q$ be a non-zero element of a field
$F$ of characteristic 0. We denote by $H_n = H_n(q,F) $ the Iwahori-Hecke algebra of
type
$A_{n-1}$, given by generators $T_i$, $i=1,2,\ ...\ n-1$ and relations
$T_iT_{i+1}T_i=T_{i+1}T_iT_{i+1}$,  $T_iT_j=T_jT_i$ if
$|i-j|>1$, and $T_i^2 = (q-1)T_i + q$.  
It is well-known that $H_n$ has a  $K$-basis 
$T_w,\ w\in S_n$, on which the generators act by
$$
T_iT_w=  
\begin{cases}
T_{s_iw} & \text{ if $\ell(s_iw)>\ell(w)$} \\ 
                   (q-1)T_{w}+qT_{s_iw} & \text{if $\ell(s_iw)<\ell(w)$}.\\
\end{cases}
$$
In fact, for any minimal expression of $w$ as a product of simple
reflections $s_i$, $T_w$ is the product of the corresponding $T_i$.

Now let $K$ be a field of characteristic 0.
In the following, we shall consider the Hecke algebras 
$H_n(x, K(x))$ over the fields $K(x)$ with parameter $q = x$, and 
the Hecke algebra
$H_n(q, K)$ over  $K$ with parameter
$q\in K^\times$. Observe that $H_n(x, K(x))$ satisfies the conditions for
the algebra $A$ in Section 0; to be consistent with the notation of Section 0, 
we shall  
denote $H_n(x, K(x))$ by $H_n$ and the specialization 
$H_n(q, K)$  by $H_n(q)$.
We sometimes refer to the algebra $H_n$ over $K(x)$ as the
{\em generic} Hecke algebra.

We recall the following definitions from Section 0: 
An element $a = \sum_w a_w(x) T_w$ of the generic Hecke algebra is said to be
{\em evaluable} at $q \in K$ if the rational functions $a_w(x)$ have no poles
at $x = q$.  A matrix representation $\varphi : H_n(x) \rightarrow \Mat_n(K(x))$
is called evaluable at $q$ if the matrices $\varphi(T_w)$ have no poles at $x = q$.
Finally, an $H_n$-module $W$ is said to be evaluable at $q$ if $W$ has a 
$K(x)$-basis such that the matrix representation with respect to this basis is
evaluable at $q$.

\v\ni
\subsection{Some induced modules}
Recall that a {\em composition} $\la$ of $n$
is  a finite  sequence $(\la_1, \la_2, \dots)$ of non-negative integers, not
necessarily decreasing,  such that $\sum_i \la_i = n$.
A {\em partition} is a composition with weakly decreasing parts.
We identify partitions with Young diagrams.

Let $\la$ be a composition of $n$ with at most $k$ non-zero parts.
Define the Young subgroup of $S_n$,
$Y_{\la}=S_{\la_1}\times S_{\la_2}\times\ ...\ \times 
S_{\la_k}\subset S_n$.
It is well-known that we can also define a length function on the
left cosets of $Y_{\la}$ by
$\ell(L_w)={\rm min}\{ \ell(v),\ v\in wY_{\la}\}$, and that
one can define a representation of $H_n(q, F)$ on the 
module $W^{\la}$, which has a basis labelled by the left cosets, by
$$
T_iL_w=  
\begin{cases}
L_{s_iw} & \text{if $\ell(L_{s_iw})>\ell(L_w)$}\\
         q L_w &\text{ if $L_{s_iw}=L_w$}\\
         (q-1)L_{w}+ q L_{s_iw} & \text{ if $\ell(L_{s_iw})<\ell(L_w)$}.\\
\end{cases}
$$
In fact, consider the subalgebra $H(\la) \subseteq H_n(q, F)$ generated by
the $T_i$ such that $s_i \in Y_\la$. $H(\la)$ acts on $F$ by
$T_i \xi = q \xi$ for all generators $T_i$ of $H(\la)$.  Then the
 module $W^\la$ is
isomorphic to
$H_n(q, F) \otimes_{H(\la)} F$, with $L_w$ corresponding to $T_w \otimes 1$.

If $\la$ and $\mu$ are composition which differ only in the order of their parts,
then the modules $W^\la$ and $W^\mu$ are isomorphic.

\v\ni
\subsection{Orthogonal representations.} Consider $H_n(q, F)$, where
$q$ is not a proper root of unity; this includes the case $q = 1$ and
also the case $F = K(x)$ and $q = x$.

It is well-known that $H_n(q, F)$ is semisimple,  its simple modules
 are labelled by
Young diagrams, and and the dimension of the simple module labelled by $\la$
is equal to the number of standard tableaux of shape $\la$.

We recall here a specific constructin of the simple $H_n(q, F)$ modules,
namely a Hecke algebra version of Young's orthogonal representations 
of the symmetric group ([H], [W1]).
Our conventions regarding Young diagrams and tableaux are those
of \cite{Macdonald}; in particular, the Young diagram $(\la_1, \la_2, \dots, \la_k)$
is a left justified array of boxes with $\la_1$ boxes in the first (top) row, 
$\la_2$ in the second, and so forth.

Let 
$$e_i=q\ident-T_i,$$
so that $e_i$ is an essential idempotent with $e_i^2 = (q+1) e_i$.
Obviously, a representation of $H_n(q, F)$ is completely determined
as soon as we know the matrices for the $e_i$'s.

For a Young diagram $\la$ of size $n$, let $V_\la$ be the vector space with  basis 
labelled by all standard tableaux of shape $\la$. 
For a standard 
tableau $t$, and for $1 \le i < j \le n$, define the quantity $d(t;i,j)$ to
be the integer, whose absolute value is one less than the length
of the hook from $i$ to $j$, and whose sign is negative
if $j$ is northeast (i.e. to the right or above) of $i$,
and positive if it is southwest of $i$. We abbreviate
$d(t;i,i+1)$ by $d(t;i)$. Moreover,  define the rational functions
$$a_d(x)=(1-x^{d+1})/(1-x^d) = {\frac{1 + x + 
\cdots+ x^d}{1+x+\cdots x^{d-1}}},\quad{\rm for}
\  d\in \N\backslash \{ 0\}.$$
A representation of $H_n(q, F)$ on $V_\la$ is determined 
by 
$$e_it=a_dt+\sqrt{a_da_{-d}}\ s_it,$$
where $d=d(t;i)$, $a_d = a_d(q)$,  and $s_i t$ is the tableau obtained by
interchanging $i$ with $i+1$. See [W1] for more details;
observe the $e_i$ here is equal to $(1+q)$ times the $e_i$ in [W1].
We will denote the representation of $H_n(q, F)$ on $V_\la$ by $\pi_\la$.

\begin{remark}
We adjoin to our field  square roots of the quantities
$a_d a_{-d}$ for $d \in \N \setminus \{0\}$.  It would also be 
possible to work with less symmetric version of these representations which avoids
the introduction of 
square roots,  namely
$$e_it=a_d t+ a_{-d}\  s_it,$$
if $a_d \ne 0$ and $e_i t = 0$ if $a_d = 0$.
\end{remark}

Note that if $i$ and $i+1$ are in the same row of $t$, then 
$e_i t = 0$, and if $i$ and $i+1$ are in the same column, then
$e_i t = (1+ q) t$.  In all other cases, both $t$ and $s_i t$ are 
standard, and the restriction of $e_i$ to the two dimensional space spanned by
$t$ and $s_i t$ has matrix
$$
\begin{pmatrix} a_d &  \sqrt{a_da_{-d}}\\
\sqrt{a_da_{-d}} & a_{-d} \\
\end{pmatrix}
$$

One can show that 
the  $V_\la$'s constitute 
complete set of inequivalent simple modules of the
semisimple algebra $H_n(q, F)$.

It is easy to obtain the restriction
rule for these modules, i.e. the way in which $V_\la$ decomposes into simple
$H_m(q, F)$-modules for
$m<n$.  One shows by induction on
$n-m$
that the multiplicity of the simple $H_{m}(q, F)$
module $V_\mu$ in $V_\la$ is equal to the number of standard skew tableaux of
shape $\la\backslash\mu$. (A standard skew tableau is a filling of the shape
$\la \setminus \mu$ with the numbers $1,2, \dots$ so that the entries of each
row and column are strictly increasing.)
\v\ni
\noindent 
\subsection{Murphy elements.} 
We shall need a $q$-version of some symmetric group
elements due to Murphy, see [M]. The
$q$-version is actually simpler than the version for the symmetric group
(see, for example, [LR]).
Let $\Delta_n^2=(T_1T_2\ ...\ T_{n-1})^n$ in $H_n(q, F)$, and
let 
$$M_i=\Delta_{i-1}^{-2}\Delta_i^2=T_{i-1}T_{i-2}\ ...\ T_2T_1^2
T_2\ ...\ T_{i-1}.$$
It is well-known that
$\Delta_n^2$ is in the center of $H_n(q, F)$; it follows easily from
this that
the elements $M_i,\ i=2,3,\ ...\ n$ are mutually commuting.
It is also well-known that the element
$\Delta_n^2$ acts by the scalar 
$\tilde\alpha_\la=q^{n(n-1)-\sum_{i<j}(\la_i+1)\la_j}$ on $S^\la$
(see e.g.  [W2, Lemma 3.2.1]).
We deduce from this that
$$M_it=q^{i-1+c_i(t) -r_i(t)}t;\eqno(*)$$
here $r_i(t)$ is the row  of $t$ in which $i$ lies 
and $c_i(t)$ the column. 
Indeed, due to the restriction rule, it suffices 
to show this for $i=n$.
In this case, $M_n$ acts on $t$ by the scalar
$\tilde\alpha_\la\tilde\alpha_{\la'}^{-1}$, where $\la'$ is $\la$
without the box containing $n$. The claim follows from an easy
computation.
\v
\ni
\subsection{Path idempotents.} We work in the generic
 Hecke algebra $H_n = H_n(x, K(x))$.
Since $H_n$ acts faithfully on $V = \bigoplus_\la V_\la$, it follows 
that there exist 
elements $p_t\in H_n$, indexed by standard tableaux of size $n$,  which are completely determined by
$p_t{\tilde s}=\delta_{t,s} t$. 
We call these elements {\em path idempotents}, because of an identification of
standard tableaux with certain paths, which we will describe later.
Recursive formulas
for these idempotents, in terms of the generators $T_i$
have been derived  in [W1, Cor 2.3].
It will be more convenient, however, to  express them in terms of Murphy
elements.

We have seen in section 1.4 that $t$ is an eigenvector of
$M_i$ with eigenvalue $\alpha_{t,i}(x)=x^{i-1+c_i(t) -r_i(t)}$.
Let $P_{t,i}$ be the eigenprojection of $M_i$ corresponding to
the eigenvalue $\alpha_{t,i}$. (Note that since the spectrum of
$M_i$ is contained in $K(x)$, it is not necessary to pass to the
algebra defined over some extension field in order to define the 
eigenprojections.)    Then we have:

\begin{lemma}
\label{ Lemma 1.5.1} $p_t=\prod_{i=2}^n P_{t,i}$
\end{lemma}
\v
\ni
\begin{proof}
It is easy to prove by induction on $i$ that
2 tableaux $t$ and $s$ are equal if and only if 
$c_i(t) -r_i(t)=c_i(s) -r_i(s)$ for $i=1,2,\ ...\ n$.
As $x$ is not a root of unity, there will be at least 
one $M_i$ which acts via distinct eigenvalues on the
the eigenvectors $t\neq s$. The result follows because the
idempotents $P_{t,i}$ are mutually commuting. \end{proof}
\v
Obviously, the Murphy elements $M_i$ are evaluable for all
$q\in K$.  If $q \in K^\times$ is not a proper root of unity, then
all path idempotents are evaluable at $q$.  This is no longer true if
$q$ is a proper root of unity.

\begin{definition} \label{Definition 1.5.2}
Fix an integer  $l > 1$.
We will say that the tableaux $t$ and $s$
are ($l$-) {\em equivalent} if 
$c_i(t) -r_i(t)\equiv c_i(s) -r_i(s) \ \mod l$ for $i=1,2,\ ...\ n$.
The {\em orbit} $[t]$ of a standard tableau $t$ consists 
of the collection of all
standard tableaux $s$ which are equivalent to $t$. We define 
the {\em orbit path idempotent} $p_{[t]} \in H_n(x, K(x))$ by
$p_{[t]}=\sum_{s\equiv t} p_s$.
\end{definition}

\v\ni
\begin{lemma}
\label{Lemma 1.5.2}
Let $q$ be a primitive $l$-th root of unity in $K$.
 The orbit path idempotent $p_{[t]}$ is evaluable
at $q$.
\end{lemma}

\begin{proof} Let $P_{[t],i}=\sum_s P_{s,i}$, where the summation goes over
all tableaux $s$ for which \break
$q^{i-1+c_i(s) -r_i(s)}=  q^{i-1+c_i(t) -r_i(t)}$.
By Proposition \ref{Proposition 0.2.3}, this idempotent is evaluable and coincides
with the eigenprojection of $M_i(q)$ for the eigenvalue
$q^{n-1+c_i(t) -r_i(t)}$.

Hence also $\prod_{i=2}^n P_{[t],i}$ is evaluable. It follows from the
definitions that it is an idempotent which projects on the span of
tableaux $s$ for which $c_i(t) -r_i(t)\equiv c_i(s) -r_i(s) \  
\mod l$ for $i=1,2,\ ...\ n$. Hence it coincides with $p_{[t]}$, which
therefore is evaluable. 
\end{proof}

\begin{lemma}  \label{Lemma 1.5.4}For any tableau  $t$ of size $n$, the path
idempotent
$p_t
\in H_n(x, K(x))$ is evaluable at $x = 1$, and the evaluation
$p_t(1)$ is the path idempotent in $KS_n$.
\end{lemma}

\begin{proof} We show this by induction on $n$, using
the recursive formulas
for the $p_t$ from [W1], Corollary 2.3.
(Since the derivation is brief, we rederive the formulas here.)
The case $n = 1$ is trivial, so suppose $n >1$ and 
fix a tableau $t$ of size $n$.
For any standard tableau $s$ of size $n$, let $s'$ be the tableau obtained
by removing the box containing the entry $n$.
We assume inductively that $p_{t'}$ is evaluable at $1$,
and that $p_t'(1)$ is the path idempotent in $KS_n$.

For standard tableaux $s$ of size $n$,
it follows from the definition of the path idempotents that
$$
p_{t'} s = 0 \quad \text{if } s' \ne t'.
$$
Moreover, it follows from 
the definition of the 
orthogonal representation that 
$$
p_{t'} e_n p_{t'}  s = 
a_{d(s)} s  \quad \text{if }  s' = t',
$$
where $d(s) = d(s; n-1)$. Hence
$$
p_t = p_{t'}  (\prod_{s} \frac{e_n - a_{d(s)}}{ a_{d(t)}- a_{d(s)}} ) p_{t'},
$$
where the product is over all $s$ such that $s \ne t$ and $s' = t'$.

We have $\displaystyle a_{d(s)}(1) =({d(s)+1})/{d(s)}$. Since $d(s) \ne d(t)$
when $s \ne t$, it follows that the product is evaluable at $1$, and therefore,
$p_t$ is evaluable at $1$.  The path idempotents in the symmetric group
algebra are given by the identical resursive formulas, with the
$a_d$ replaced by the rational numbers $(d+1)/d$; so it follows from
the induction hypothesis that $p_t(1)$ is the path idempotent in 
the symmetric group algebra for the tableau $t$.
\end{proof}

\v
\subsection{Specht modules.} The 
$H_n(x, K(x))$-modules $V_\la$ 
are not in general evaluable  at $q$ when
$q$ is a proper root
of unity. So-called Specht modules
$S^\la$ have been defined in [DJ] 
which are equivalent to $V_\la$ (as $H_n(x, K(x))$ modules) and which
are evaluable also at roots of unity. We briefly recall the
construction and properties of these modules; see [DJ] for more
details. 

The construction of the modules $S^\la$ can be 
carried out for $H_n(q, F)$ for any $q$ and $F$.
It turns out that for any $q \in K^\times$,
the ``restriction" of the
evaluable $H_n(x, K(x))$ module $S^\la$ to $H_n(q)$ coincides with module $S^\la$
constructed over $H_n(q)$. It is safe, therefore, to use the same notation
for all of these modules.

Let $Y_{\la}=S_{\la_1}\times\ ...\ 
\times S_{\la_r}\subset S_n$.
We define the symmetrization and antisymmetrization operators
$$\Sym_\la=\sum_{w\in Y_{\la}} T_w\quad {\rm and}\quad
A_\la = \sum_{w\in Y_{\la}} (-q)^{n(n-1)/2-\ell(w)}T_w.$$
Let $S^{\la}$ be the submodule of the $H_n(q, F)$-module 
$W^{\la}$ (see 1.2) generated
by $A_{\la' }W^{\la}$,
where $\la'$ is the diagram obtained from $\la$ by interchanging
rows with columns. Hence $S^\la \cong H\,A_{\la'}\,H\,\Sym_\la$,
where $H$ denote $H_n(q, F)$.

\begin{remark}
Let $w_\la$ be the permutation which maps the standard tableau of
shape $\la$ filled column by column (i.e. the first column 
is filled first, then the second etc) to the standard tableau
of shape $\la$ obtained by filling $\la$ row by row.
(This $w_\la$ is the inverse of the $w_\la$ in [DJ].)
By [DJ, Lemma 4.1], $S^\la \cong H_nA_{\la '}T_{w_\la}\Sym_\la$;
the latter is the Specht module in the definition of [DJ, Section 4.1].
\end{remark}

\v\ni
\begin{theorem}
\label{Theorem 1.6.1} {\rm (Dipper-James)} Write $H_n = H_n(x, K(x))$.
\v\ni
\begin{enumerate}
\item
There exists an explicit
basis of  the $H_n$-module 
$H_n\,A_{\la'}\,H_n\,\Sym_\la  \subseteq H_n$ which is evaluable for all $q \in
K^\times$,  and such that the basis elements evaluated at $q$ remain
linearly independent over $K$ for all $q\in K^\times$. 
\item
The action of the generators
$T_i$ on $S^\la \cong H_n\,A_{\la'}\,H_n\,\Sym_\la $ is given by evaluable matrices
with respect to this basis. 
\item
The restriction $S^\la_q$ of the Specht module constructed over
$H_n$ is identical with the Specht module constructed
over $H_n(q)$.
\item
The $H_n$ module $S^\la$ is equivalent over $K(x)$ to the $H_n$-module
$V_\la$ (defined in 1.3).  
For $q$ not a proper root of unity, the $H_n(q)$ module
$S^\la$  is equivalent over $K$ to the $H_n(q)$-module  $V_\la$.
\end{enumerate}
\end{theorem}
\v
\ni
\begin{proof} This theorem is   [DJ, Theorem 5.6].
The statements (b) and (d) follow easily 
 by carrying out  construction of [DJ]  over the domain $K[x]$
of polynomials in $x$. \end{proof}

In the following, we use two orderings on the set of Young diagrams
of a given size.  The dominance order is defined by
$\mu \trianglerighteq \la$
if $\sum_{i = 1}^r \mu_i \ge \sum_{i = 1}^r \la_i$ for all $r$.
This is a partial order. 
The lexicographic order is defined by $\mu > \la$ if the first non-zero
difference $\mu_i - \la_i$ is positive.  This is a total order.
One has $\mu \trianglerighteq \la \Longrightarrow \mu \ge \la$.
\v\ni
\begin{theorem}
\label{Theorem 1.6.2} {\rm(Dipper-James)} \quad
Let $q \in K$ be a primitive $l$-th root of unity. 
Consider the Hecke algebra $H_n(q)$ over $K$.
\v\ni
\begin{enumerate}
\item
Let $(\ ,\ )$ be the bilinear 
form on the $H_n(q)$-module $W^{\la}$ with orthonormal
basis $(L_w)$. Then the module $D^\la=S^\la/(S^\la \cap (S^\la)^\perp)$
is either 0 or it is simple.
\item
 A diagram $\mu$ is called $l$-regular if it has at most $l-1$
rows of equal length. The module $D^\mu$ is nonzero if and only if
$\mu$ is $l$-regular {\rm ([DJ, Theorem 6.3 and 6.8])}, 
and the set $\{D^\mu :\mu\  {\rm  is  }\ l{\rm -regular}\}$
is a complete set of mutually nonisomorphic simple
of $H_n(q)$ modules {\rm ([DJ, Theorem 7.6])}.
\item
 Let $0\subset V_1\subset V_2\subset
\ ...\ V_r=W^\la$ be a composition series of $W^\la$ such that the factors
$V_{i+1}/V_i$ are simple $H_n(q)$ modules. Then any simple factor 
has to be isomorphic to $D^\mu$ for some $\mu \trianglerighteq \la$
(in dominance order).
If $\la$ is $l$-regular,   there is
exactly one factor in such a series which is isomorphic to $D^\la$
{\rm ([DJ, Theorem 7.6])}. In particular, if $\la$ is $l$-regular,
$D^\la$ appears with multiplicity 1 in any composition series of $S^\la$.
\end{enumerate}
\end{theorem}

\v
\noindent 
\subsection{Decomposition numbers.} 
Recall the notions of the rank vector of an idempotent from subsection 0.5.
We say that $\la$ is the $highest\ component$ of the idempotent $p$
if $\la$ is the highest diagram (in lexicographic order) for which 
$r(p)_\la\neq 0$.
Let $q$ be a primitive $l$-th  root of unity in $K$. 
If $p$ is an idempotent of $H_n(q)$,
its $q$-rank ${\bf r}_q(p)$ in the maximum semisimple
quotient $\overline{H_n(q)}$ of $H_n(q)$ is given by
a vector $(r_q(p)_\mu)$, where $\mu$ runs through all
$l$-regular diagrams. 

As in section 0.6, we define numbers $d_{\la\mu}$
by $\dlm = {\bf r}(p_\mu)_\la$, where $p_\mu$ is an evaluable
idempotent in $H_n$ whose image in $\overline{H_n(q)}$
is a minimal
idempotent in the component $\overline{H_n(q)_\mu}$.
We can now reformulate the results of  Section 1.6 in the following 
way:

\v\noindent
\begin{proposition}
\label{Proposition 1.7.1} 
\quad Let $q$ be a primitive $l$-th root of unity in $K$.
 Let $p$ be an idempotent in
$H_n$ which is evaluable at $q$,  and let
 $\mu$ be the highest component of $p$. Then we
have:
\smallskip\noindent
\begin{enumerate}
\item
The highest nonzero component of ${\bf r}_q(p)$ is also $\mu$.
In particular, $\mu$ has to be $l$-regular. 
\item
$\Tr_{S^\mu}(p(q))=\Tr_{D^\mu}(p(q))$.
\smallskip\noindent
\item
$d_{\la\mu}\leq {\bf r}(p)_\la/{\bf r}(p)_\mu$ for all $\la$.
\item
$d_{\la\mu}$ is the number of factors isomorphic to $D^\mu$
in any composition series of $S^\la$ as in Section 1.6.
\end{enumerate}
\end{proposition}

\v
\ni
\begin{proof} By assumption, $p$ acts as 0 on $S^\nu$ for any $\nu>\mu$,
and hence also its evaluation $p(q)$; in particular $p(q)$ acts
as 0 on $D^\nu$. As $p$ and $p(q)$ act nonzero
on $S^\mu$, $D^\mu$ has to be nonzero, and $p(q)$ has to act
on it as a nonzero endomorphism, by Theorem \ref{Theorem 1.6.2}(c). This shows (a).
Statement (b) follows from this and \ref{Theorem 1.6.2}(c).

To prove (c), let $\overline{p(q)}$ denote the image of 
$p(q)$ in $\overline{H_n(q)}$, and decompose $\overline{p(q)}$ as
$\overline{p(q)} = \overline{p(q)}_\mu + \sum_{\nu <
\mu}\overline{p(q)}_\nu$. Let $p_\mu$ be a lifting of $\overline{p(q)}_\mu$
to an evaluable idempotent in $H_n$ satisfying 
$p_\mu = p p_\mu p$.  Then we have 
$$
\Tr_{D^\mu}(p_\mu(q)) = \Tr_{D^\mu}(p(q)) = \Tr_{S^\mu}(p(q)) = {\bf r}(p)_\mu,
$$
using (b).  Also
$$
\Tr_{S^\la}(p_\mu) \le \Tr_{S^\la}(p) = {\bf r}(p)_\la.
$$
Therefore 
$$d_{\la \mu} = {\frac{\Tr_{S^\la}(p_\mu)}{\Tr_{D^\mu}(p_\mu(q))}}
\le {\bf r}(p)_\la/{\bf r}(p)_\mu. $$
Statement (d)  follows from Proposition \ref{Proposition 0.6.1} and Theorem 
\ref{Theorem 1.6.2}. \end{proof}

\v\ni
\begin{remark}
By the last proposition, in order to get upper bounds for the 
coefficients
$d_{\la\mu}$,
it suffices to construct evaluable idempotents $p$ in $H_n$ for which
${\bf r}(p)_\nu=0$ for all $\nu > \mu$   This will be our strategy
in sections 3 and 4.
\end{remark}
\v

\section{The $k$-row quotient.}
\subsection{Definition of the $k$-row quotient}

\begin{definition} Fix  natural numbers $k \le n$.  
Let $W(n,k) = W(n,k,q,F)$ be the
direct sum of all $H_n(q,F)$-modules $W^\la$,
  for  Young diagrams $\la$  of size
$n$ with at most
$k$-rows. The $k$-row quotient
of the Hecke algebra $H_n(q, F)$ is the quotient by the kernel of
the module $W(n,k)$.
\end{definition} 

We remark that the $k$-row quotient can also be described in terms of
a representation of the Hecke algebra on tensor space which
has a natural connection with quantum groups of type $A$.

 Let $V$ be a $k$-dimensional vector space over $F$, 
with basis $\vareps_1, \vareps_2, \ ...
\ \vareps_k$. Then $V^{\otimes n}$ has a natural basis
of the form $\vareps_{i_1}\otimes\ ...\ \otimes \vareps_{i_n}$, where
$1\leq i_j\leq k$
for $j=1,2,\ ...\ n$. 
Define a matrix $R$ in $\End(V\otimes V)$ by
$$R(e_i\otimes \vareps_j)=
\begin{cases} \vareps_j\otimes \vareps_i & \text{if $i<j$}\\
                  q\, \vareps_i\otimes \vareps_i &\text{if $i=j$}\\
                   (q-1)\,\vareps_i\otimes \vareps_j+ q\,\vareps_j\otimes 
				\vareps_i &\text{ if $i>j$}\\
\end{cases}
$$
Using these matrices $R$,  define operators $R_i$ in 
End $(V^{\otimes n})$ by
$$R_i = 1\otimes 1\otimes\ ...\ 1\otimes R\otimes 1\otimes\ ...\ \otimes 1,$$
which acts by the matrix $R$ on the $i$-th and $(i+1)^{st}$ factor of
$V^{\otimes n}$. This is the famous representation
of the Hecke algebra $H_n(q, F)$ 
 discovered by Jimbo [Ji]. The fact
that it is indeed a representation is a consequence of the following observations.

  For any composition $\la$ of
$n$  with at most $k$ non-zero parts define
the vector   
$$ \vareps^{\, \otimes \la} = \vareps_1^{\otimes \la_1}\otimes\ ...
\ \otimes \vareps_k^{\otimes \la_k},$$
where $\vareps_i^{\otimes m}= \vareps_i\otimes\ ...\ \otimes \vareps_i$, 
($m$ times), and
$$V^{\otimes \la} = {\rm span}\ \{ w(\vareps^{\,\otimes \la}),
\ w\in S_n\} ,$$
where the symmetric group $S_n$ acts by permuting the factors in
the tensor product.
It is clear that $V^{\otimes n}$ is the direct sum, over all compositions $\la$
of
$n$ with at most $k$ non-zero parts, of $V^{\otimes \la}$, and that
$V^{\otimes \la}$ is invariant under the operators $R_i$.
We define
$$\ell(\vareps_{i_1}\otimes\ ...\ \otimes \vareps_{i_n}) = 
\# \{ (i_j ,i_k) : j<k\ {\rm and}\  i_j > i_k\ \ \} ,$$
i.e. the number of inversions of indices.
As the stabilizer of $\vareps^{\otimes \la}$ is equal to
$Y_{\la}$ (see section 1.2), the map $\iota :
L_w \mapsto w(\vareps^{\otimes \la})$ induces an isomorphism
between the vector spaces $W^{\la}$ and $V^{\otimes \la}$, preserving
the length functions defined on the basis vectors
of these two vector spaces.
Using this observation, it is easy to check that
$\iota$ intertwines the action of $T_i$ on $W^{\la}$,
and of $R_i$ on $V^{\otimes \la}$. Thus $V^{\otimes \la}$ is an
$H_n(q, F)$ module, isomorphic to $W^\la$.

Finally, $T_i \mapsto R_i$ defines a representation of $H_n(q, F)$ on
$V^{\otimes n}$, and 
$$
V^{\otimes n} = \bigoplus V^{\otimes \la} \cong \bigoplus W^ \la,
$$
where the sum is over all compositions $\la$ of $n$ with at most $k$ non-zero
parts.  

Given a composition $\la$ of $n$ let $\mu$ be the partion whose parts are the
same as those of $\la$, but in  decreasing order.  Then $W^\la \cong W^\mu$.
Thus
$$
V^{\otimes n} \cong \bigoplus m_\la W^ \la,
$$
where now the sum is over partitions of $n$ with at most $k$ parts, 
and the $m_\la$ are certain multiplicities.  Thus the quotient of
$H_n(q, F)$ which acts faithfully on $V^{\otimes n}$ is the same as the 
$k$-row quotient defined above.

We  denote by $\Phi_k$ the   matrix representation of 
$H_n(x, K(x))$ on  $(K(x)^k)^{\otimes n}$ with respect to the
natural basis described above.   This
representation is  evaluable at any non-zero $q \in K$, and the
``restriction" to 
$H_n(q)$ is again the representation on  $(K^k)^{\otimes n}$ with respect
to the natural basis. 

Denote by $V(n,k)$ the direct sum of the $H_n(x, K(x))$ modules
$V_\la$ for Young diagrams $\la$ of size $n$ with at most $k$ rows,and
by $\Pi_k$ the matrix representation of $H_n(x, K(x))$ on $V(n,k)$ with
respect to the basis of standard tableaux.  Since the simple direct summands
of $(K(x)^k)^{\otimes n}$ are those labelled by Young diagrams with no more than $k$
rows,   $\Pi_k(H_n(x, K(x)))$ is isomorphic to the $k$-row quotient of 
$H_n(x, K(x))$.

\v
\ni
\subsection{Affine Weyl group, $l$-cores, and blocks.}

Fix an integer $k\ge 2$, and $W$ be the affine reflection group of type
$A_{k-1}^{(1)}$. It is isomorphic to the semidirect product $\Z^{k-1}\times S_k$.
Fix $l\in\N$. Then we have a faithful action of $W$ on $\R^k$,
with $S_k$ acting via permutation of the coordinates, and with 
$$\tau_i(y_1,\ ...,\ y_k)= (y_1,\ ...,\ y_i-l,y_{i+1}+l,\ ...\ y_k),$$
where $\tau_i$ is the $i$-th generator of the translation group
$\Z^{k-1}\subset W$. Observe that any affine reflection is
a reflection in a hyperplane given by an equation of the form
$y_i-y_j=ml$ for some $i,j$ with $1\leq i,j\leq k$ and with $m\in\Z$.

Define $\rho = (k-1,k-2,\ ...,\ 1,0)$. Another important action 
(the ``dot action") of
$W$ on $\R^k$ is given by $$w\cdot v = w(v + \rho) - \rho.$$
We next review the relation between the orbits of the dot action of $W$ 
and the notion of the $l$-core of a diagram.

We recall  the
notions  of a rim hook and $l$-core, for example from [JK].
For $(a,b)$ a node of a Young diagram $\la$, the corresponding
{\em hook} is the portion of row $a$ to the right of $(a,b)$ together
with the portion of  column $b$ below $(a,b)$,   including the cell $(a,b)$;
the length of the hook is $h_{(a,b)} = \la_a -a + \la'_b - b + 1$.
The corresponding {\em rim hook} is the portion of the rim of $\la$
between $(\la'_b,b)$ and $(a, \la_a)$; the rim hook is a connected
skew shape of size $h_{(a,b)}$.
 A diagram is called an $l$-core if it has no rim hook of length $l$. 
Every Young diagram contains a unique $l$-core, which is obtained by 
removing successive rim hooks of length $l$.

Let $D$ be 
the set of points $y \in \R^k$ such that $y_1>y_2>y_3>\ ...\ >y_k$.
Let $D^+$ be the set of points $y \in D$ with all components
nonnegative.
We obtain a map from the set of Young diagrams with $k$ rows at the
most onto $D^+\cap\Z^k$, given by 
$\la\mapsto \la+\rho$,
where on the right hand side we identify $\la$ with the vector
$(\la_i)_i$.
Let $\la$ be a diagram  with at most $k$ rows. The point
$y = \la + \rho \in D^+$ is determined by its set of components
$\{y_i\}$.  It is not difficult to see that the operation of
removing a rim hook of length $l$ from $\la$ corresponds exactly
to the operation of reducing one element of $\{y_i\}$ by $l$,
see [JK], Lemma 2.7.13.

The following result must be well known.

\vbox{
\begin{lemma}
\label{Lemma 2.4.1} Let $\la$ and $\mu$ be diagrams of the same size with
at most $k$ rows.
Let $y = \la + \rho$, $z = \mu + \rho$.  The following are equivalent:

{\rm(1)}\quad $\la$ and $\mu$ have the same $l$-core.

{\rm(2)}\quad For $0 \le r \le l-1$, \ $|\{i : y_i \equiv r\ \  \mod\ \  l\}|
= |\{i : z_i \equiv r\ \  \mod\ \  l\}|$.

{\rm (3)} \quad $\mu$ is in the  orbit of $\la$ under the dot action of $W$.
\end{lemma}
}
\v\ni
\begin{proof}
The equivalence of statements (1) and (2) follows from
the observation made just above.  Furthermore, (3) $\Longrightarrow$ (2)
is evident since an affine reflection changes each of two
coordinates of $y$ by a multiple of $l$, so preserves the 
cardinalities of the sets $J_r = \{i : y_i \equiv r\ \  \mod\ \  l\}$.
So it remains to show (2) $\Longrightarrow$ (3).  One can define a
distance between $\la$ and $\mu$ as follows:  Let $y_{i,r}$ be
the elements of $J_r$ in increasing order, and similarly let
$z_{i,r}$ be the elements of $J'_r = \{i : z_i \equiv r\ \  \mod\ \  l\}$
in increasing order. Define $|\la - \mu| = \sum_r \sum_i |y_{i,r} -z_{i,r}|$.
It suffices to
show that if $\la \ne \mu$, then there is an element $w \in W$ such
that $w\cdot\la$ is a Young diagram and $|w\cdot \la - \mu| < |\la - \mu|$.

One has $\sum_r \sum_i y_{i,r} = \sum_r \sum_i z_{i,r}$.
Suppose that there is some $r$ such that 
$ \sum_i y_{i,r} \ne  \sum_i z_{i,r}$.  Then one can check that there
is an element $w \in W$ and $s \ne t$ such that
$w\cdot\la$ is a Young diagram,  $w$ increases one element of $J_s$ and
decreases one element of $J_t$, and $|w\cdot\la - \mu| < |\la - \mu|$.

Suppose one the other hand that $ \sum_i y_{i,r} =  \sum_i z_{i,r}$
for all $r$.  Choose $r$ such that $J_r \ne J'_r$.  Then one can
show that there is an element $w \in W$ such that $w\cdot\la$ is a Young diagram,
$w$ increases one element of $J_r$ and decreases another, and
$|w\cdot\la - \mu| < |\la - \mu|$. 
\end{proof}

\v

The blocks of the Hecke algebras of type A at a root of unity
are parametrized by  $l$-cores:  
\v\ni
\begin{theorem}
\label{Theorem 2.4.2} {\rm ( [DJ2], Theorem 4.1.3)} \quad  Let $q$ be a 
primitive $l$-th root of unity in
$K$.
\v\ni
\begin{enumerate}
\item
If $\la$ and $\mu$ are two $l$-regular diagrams of size $n$, then
$D^\la$ and $D^{\mu}$ belong to the same block of $H_n(q)$ if and only if
$\la$ and $\mu$ have the same $l$-core.
\item
If $D^\mu$ is a composition factor of $S^\la$, then
$\mu \trianglerighteq \la$ and $\mu$ and $\la$
have the same $l$-core.
\end{enumerate}
\end{theorem} 

\v\ni
\subsection{Paths and path equivalence.}
If $t$ is a standard tableau, 
we define $t(i)$ to be the Young diagram consisting of the
boxes containing the numbers 1, 2, ..., $i$.
If $t$ has $n$ boxes and at most $k$ rows, we identify  $t$ with
the piecewise affine path $t:[0,n]\to \R^k$ which takes the values $t(i)$
at $i = 0, 1, \dots, n$ and which is affine on $[i, i+1]$.
We denote by $t_\rho$ the path $t_\rho(s) = t(s) + \rho$.
 Observe that
$t_\rho([0,n])\subset D^+$.  Likewise, if $t$ is a standard skew tableau 
of shape $\la \setminus \nu$, 
we identify $t$ with a piecewise affine path from $\nu$ to $\la$, and 
we also consider the path $t_\rho(s) = t(s) + \rho$ which goes from
$\nu + \rho$ to $\la + \rho$.

Let $t$ be a standard (skew) tableau with at most $k$ rows.
Let $t_\rho(i)_a$ denote the $a$-th coordinate
of $t_\rho(i)$. Recall that    $c_j(t)$and $r_j(t)$ denote
the column and row of the box of $t$ which contains the number $j$.
Using these notations, 
$c_i(t)-r_i(t)=t_\rho(i)_{r_i}-k.$
Hence, the tableaux $t$ and $\tilde t$ are equivalent in the sense of
Definition \ref{Definition 1.5.2}  if and only if 
$$t_\rho(i)_{r_i(t)}\equiv
\tilde t_\rho(i)_{r_i(\tilde t)}\quad {\rm mod}\  l\quad {\rm for\ all\ }i.\eqno(*)$$ 
This definition of $l$-equivalence obviously extends to 
arbitrary piecewise linear paths \break $p:[0,n]\to \R^k$ for which $p(i)\in \Z^k$.
The notion of $l$-equivalence can be characterized
by defining something like an `action' of the affine reflection
group on paths,  as follows:
Assume that $t_\rho(i)$ lies in an affine hyperplane corresponding
to an affine reflection $s$ (as defined in Section 2.4).
Then we define 
$$s^{(i)}(t_\rho(u))=
\begin{cases} t_\rho(u) & \text{if $0\leq u\leq i$,}\\
         s(t_\rho(u))& \text{ if $i<u\leq n$.}\\
\end{cases}
$$
Geometrically, it means we reflect the part of the path $t_\rho$
after the point $t_\rho(i)$ in the hyperplane belonging to $s$.
The index $i$ is necessary, as $t$ might cross the  hyperplane 
of $s$ more than once. Observe that $s^{(i)}(t_\rho)$
may be a path which no longer belongs to a tableau, but it remains
equivalent to $t_\rho$. 
\v\ni
\begin{lemma}
\label{Lemma 2.5.1} Two paths $t$ 
and $\tilde t$ are $l$- equivalent  if and only
if there exists a finite
sequence of affine reflections $s_1,\ ...\ s_r$ and integers
$i_j\in [1,n]$, $j=1,\ ...,\ r$ such that 
$$\tilde t_\rho=s_1^{(i_1)}(s_2^{(i_2)}(\ ...\ s_r^{(i_r)}(t_\rho)\ ...\ )).$$
\end{lemma}
\v
\ni
\begin{proof}
Assume that $t$ and $\tilde t$ are $l$-equivalent and 
assume that  $t_\rho$ and $\tilde t_\rho$ coincide 
from 0 to $i$, but not at $i+1$. Let $a$ and $b$ be the rows in which
the box containing $i+1$ is added in $t$ resp. $\tilde t$.
Then $t_\rho(i+1)_a\equiv \tilde t_\rho(i+1)_b \  \mod l$,
which implies $t_\rho(i)_a\equiv \tilde t_\rho(i)_b=t_\rho(i)_b \  \mod l$.
Hence $t_\rho(i)$ lies on a hyperplane of the form $y_a-y_b=ml$, for some
$m$, and the paths $\tilde t_\rho$ and $s^{(i)}(t_\rho)$ coincide
up to $i+1$ if we define $s$ by
$$s(y_1,\ ...,\ y_a,\ ...\ y_b,\ ...,\ y_n)=
(y_1,\ ...,\ y_b+ml,\ ...\ y_a-ml,\ ...,\ y_n).$$
Repeat this construction for the paths $\tilde t_\rho$ and 
$s^{(i)}(t_\rho)$ until the resulting paths coincide. 
\end{proof}
\v
\subsection{Some simple Specht modules}
The following two lemmas present very simple sufficient criteria
for a Specht module to be simple.  A more general criterion is
given below in Corollary \ref{Corollary 3.3.4}.  All of these results are special
cases of results of James and Dipper ([DJ], 4.11) and James and
Mathas [JM].
\v\ni
\begin{lemma}
\label{Lemma 2.6} Assume that $\la$ is the highest diagram of its
orbit under the dot action of the affine Weyl group $W$.
Then $D^\la=S^\la$.
\end{lemma}
\v
\ni
\begin{proof}
If $\mu \trianglerighteq \la$ and $\mu$ has the same $l$-core as 
$\la$, then it follows that $\mu$ has at most $k$ rows and
furthermore $\mu$ is in the orbit of $\la$ under $W$.
By assumption, $\mu = \la$.  It follows that $D^\la$ is the only 
composition factor of $S^\la$.
\end{proof}
\v\ni
\begin{lemma}
\label{Lemma 2.7}  Let $\la = (l-1) \rho$.
Then 
$S^\la = D^\la$ and for all diagrams $\nu$ such that $\nu \ne \la$
and $|\nu| = |\la|$, $d_{\nu \la} = 0$.
\end{lemma}
\v\ni

\begin{proof}
This follows since $\la$ is its own $l$-core.
\end{proof}

\section{Big diamond elements and their matrix coefficients}

In this section we show that the $k$-row quotient of $H_n(x)$
contains a (non unital) subalgebra isomorphic to the $k$-row 
quotient of $H_m(x^l)$, when $n = ml + (l-1){\binom{k} {2}}$. 
Furthermore, when $q$ is an $l^{\rm th}$ root of unity, then the
$k$-row quotient of $H_n(q)$
contains a (non unital) subalgebra isomorphic to the $k$-row 
quotient of $KS_m$.  This is our main technical result.

As a corollary of this result, we obtain the simplicity of certain
Specht modules over $H_n(q)$, in characteristic 0 (a special
case of results of Dipper and James, and of James and Mathas.)

\v\ni
\subsection{$l$-Straight tableaux and $k$-critical diagrams.}
\begin{definition}
We say that
a skew tableau of shape $\la\backslash\nu$ 
 is {\em $l$-straight} if 
$$d(t;i,i+1)\equiv -1 \quad  \mod l$$ for $i=1,2,\ |\la\backslash\nu |$.
\end{definition}

This condition can be expressed as $$(c_{i+1}(t) - r_{i+1}(t)) -(c_i(t) - r_i(t))   
\equiv
1\ \mod l$$ for all $i$.  Geometrically, the condition says that the path
$t$ is equivalent to a path in which only one coordinate is increased, that is to
a straight path.

\begin{definition} Fix an integer $k \ge 2$.  
A {\em point} $y\in D^+\cap \Z^k$ is said to be {\em $k$-critical}, or just
critical, if $y_i-y_j$ is divisible by $l$ for all $1\leq i,j\leq k$;
it is called a {\em reduced $k$-critical point} if
it is $k$-critical and $y_k=0$.
A {\em diagram} $\la$ with at most $k$ rows is said to be {\em $k$-critical}
 if $\la + \rho$
is a $k$-critical point, that is if $\la_a - \la_b + b-a$ is divisible by
$l$ for all $a$ and $b$.  It is  said to be a {\em reduced $k$-critical
diagram} if in addition $\la_k = 0$.
\end{definition}

The smallest 
$k$-critical point 
is $l \rho$, and the smallest $k$-critical Young diagram is
$(l-1) \rho$.  The size of the smallest  $k$-critical Young diagram
is 
$n_0 = (l-1){\binom{k}{2}}$.

\begin{definition}
We call an $l$-straight skew tableau $t$ {\em $k$-special} if the initial
shape of $t$ is $k$-critical and if the length of $t$ is a multiple of $l$.
\end{definition}

If $t$ is a $k$-special  $l$-straight skew tableau with initial shape $\mu$ and
length $n$,
then for all $i$,  at most one coordinate of $t_\rho(i) - \mu_k =
t(i) + \rho - \mu_k$ is not divisible by $l$, and that coordinate is
congruent to $i \ \mod  l$, as one can see easily by induction.  Consquently
 $t_\rho(ml)$ is a $k$-critical point for all $m$ with
$ml\leq n$.

Let $T^l$ denote the set of $k$-special $l$-straight skew tableaux with
initial shape $(l-1)\rho$, and,   for
a $k$-critical diagram $\la$, let $T^l_\la$ denote the set of $t \in T^l$
which end in $\la$.

If $t \in T^l$ with length $n$,  then for all $m$ such that $ml \le n$,
the diagram 
 $t_\rho(ml)/l$ has strictly decreasing components,
so $\tilde t(m) = t_\rho(ml)/l - \rho$ is
 a Young diagram with $m$ boxes.
It is easy to check that the sequence of Young diagrams $\tilde t(m)$
defines a Young tableau of shape $t_\rho(n)/l-\rho$. This shows the
following lemma.

\v\ni
\begin{lemma}
\label{Lemma 3.1.1} Let $\la$ be a  $k$-critical diagram. There is 
a bijection  $\Psi$ 
between $T^l_\la$ and the set of all standard Young tableuax
 of shape $(\la+\rho)/l - \rho$
defined by $\Psi(t)(m) = t_\rho(ml)/l - \rho$. \qed \end{lemma}

\v
\subsection{Big diamond elements}\  Recall that $e_i= x\ident - T_i$ in the
generic Hecke algebra over $K(x)$.
A big diamond element is one of the form
$$
E(n) = (e_{n+l}e_{n+l+1}\ ...\ e_{n+2l-1})(e_{n+l-1}e_l\ ...\ e_{n+2l-2})\ ...
\ (e_{n+1}e_{n+2}\ ...\ e_{n+l}).
$$
Here $l$ is fixed,  and $n$ is
arbitrary.
For simplicity of notation we will often normalize $n$ to be zero when
discussing big diamond elements, and write $E$ for the big diamond element.
For  any finite set $A$ of natural numbers containing
no two consecutive numbers,  let $e_A$ be the product of the commuting
elements $e_j$ for $j\in A$ and $s_A$ the product of the 
simple transpositions $s_j$ for $j \in A$.
Let  $J_i = \{l-i+1, l-i+3,\dots, l+i-1\}$ for $1\le i\le l$.
Thus  $J_1 = \{l\}$, $J_2 =  \{l-1, l+1\}$, etc.
We define $e_J=e_{J_{l-1}}e_{J_{l-2}}\ ...\ e_{J_2}e_{J_1}$
and $s_J =  s_{J_{l-1}}\dots s_{J_{2}} s_{J_{1}}$.
We also write $e_{J,i} = e_{J_i} e_{J_{i-1}}\cdots e_{J_2} e_{J_1}$,
and   $s_{J,i} = s_{J_i} s_{J_{i-1}}\cdots s_{J_2} s_{J_1}$.
Then one can check that
$$
E 
= e_{J_1}e_{J_2}\dots e_{J_{l-1}}e_{J_l}e_{J_{l-1}}
\dots e_{J_2} e_{J_1}.
$$
Let $s_E = s_{J,l-1}\inv s_{J_l} s_{J,l-1}$; thus
$$s_E(i)=
\begin{cases}
 i+l & \text{ if $i\leq l$,}\\
i-l&\text{if $l<i\leq 2l$.}\\
\end{cases}
$$

For a pair of $k$-special $l$-straight skew tableaux 
$t_1$ and $t_2$ with the same initial
and final shapes, one has
\begin{equation} \label{Equation 3.2.1}
\langle Et_1, t_2\rangle
\ =\ \langle e_{J_l}e_Jt_1, e_Jt_2\rangle.
\end{equation}
The bilinear form used here is the one in which standard skew tableaux are
orthonormal.
Our goal is to obtain a formula for such matrix coefficients.

In the following, let 
$t$ be a $k$-special $l$-straight skew tableau of length $2l$. 
Let $\mu$ denote the ($k$-critical) initial shape of $t$; note that
the first $l$ cells of $t$ are added to $\mu$ in some row  $a$
and the next $l$ cells in some row $b$ (possibly $a = b$).
Write  $dl = \mu_a - \mu_b + b-a$.

For any integer $d$, we write $a_{\pm d} = a_d a_{-d}$,

\v
\noindent
\begin{lemma}
\label{Lemma 3.2.2}\quad 
For $1 \le i \le l-1$,
$$
e_{J_{i+1}} e_{J_i} s_{J,i-1} t =  {a_{\pm(dl + l -i)}}^{i/2}\
 e_{J_{i+1}} s_{J,i} t.
$$
\end{lemma}

\medskip\ni
\begin{proof} By definition of the orthogonal representation, one
has
$$
e_{J_{i+1}} e_{J_i} s_{J,i-1} t =  
\sum_{H \subseteq J_i}\alpha_H e_{J_{i+1}} s_H s_{J,i-1} t,
$$
where $\alpha_H^2$ is a rational function.  We have to show that
if $H \ne J_i$, then  $e_{J_{i+1}} s_H s_{J,i-1} t = 0$.

Observe that the skew tableau $s_{J,i-1} t$ has the digits
1 through $l-i$ in row $a$, followed by the elements
 of $J_{i}$ in increasing order, while 
in row $b$, the skew tableau   has  
$\{1 + r : r \in J_{i}\}$ in increasing order followed by 
$l+i+1, \dots, 2l-1, 2l$.  

If $H \ne J_i$, let $r$ be the least element of $J_i$ such that $r \not\in H$.
Then $s_H s_{J,i-1} t$ has the digits $r-1$ and $r$ in successive cells
in row $a$.  Since $r-1 \in J_{i+1}$, it follows that
 $e_{J_{i+1}} s_H s_{J,i-1} t = 0$.

>From the description of $s_{J,i-1} t$, one sees that for $r \in J_i$,
$$d(s_{J, i-1} t, r, r+1) = \pm(\mu_a - \mu_b + b-a + l-i) 
= \pm(dl + l-i).
$$
Thus  
$$
\langle s_{J,i} t, e_{J_i} s_{J,i-1} t\rangle = {a_{\pm(dl + l -i)}}^{i/2}. 
$$
\end{proof}

\v
\ni
\begin{lemma}
\label{Lemma 3.2.3}  For $r, s \ge 1$, 
$$
(a_{\pm r})^{s+1} (a_{\pm(r+1)})^s \cdots a_{\pm(r+s)} = 
x^{\binom{s+2} 2} a_{r-1}^{-(s+1)}\ {\frac{1-x^{r+s+1}} {1-x^r}}.
$$
\end{lemma}
\medskip\ni
\begin{proof}
One has 
$$a_r a_{r+1}\cdots a_{r+s} = {\frac{1-x^{r+s+1}}{1-x^r}},$$
and therefore

\begin{equation}\label{Equation 3.2.2}
\begin{split}
(a_r)&(a_r a_{r+1})\cdots (a_r  a_{r+1}\cdots a_{r+s}) = \\
 &(1-x^r)^{-(s+1)}
(1-x^{r+1})(1-x^{r+2})\cdots(1-x^{r+s+1}).
\end{split}
\end{equation}
Similarly, 
$$
a_{-r} a_{-r-1}\cdots a_{-r-s} = x^{s+1}
{\frac{1-x^{r-1}}{1-x^{r+s}}},
$$
so

\begin{equation}\label{Equation 3.2.3}
\begin{split}
(a_{-r})&(a_{-r} a_{-r-1})\cdots (a_{-r}  a_{-r-1}\cdots a_{-r-s}) = \cr
&x^{\binom{s+2} 2} (1-x^{r-1})^{s+1}(1-x^r)\inv(1-x^{r+1})\inv\cdots (1-x^{r+s})\inv. 
\end{split}
\end{equation}
The result follows from multiplying the expressions in Equations
\ref{Equation 3.2.2} and \ref{Equation 3.2.3}.
\end{proof}

\v\ni
\begin{lemma}
\label{Lemma 3.2.4} 
 $$
e_{J_l} e_{J, l-1} t = \left(x^{\binom{l}{2}} 
a_{dl}^{-l+1}\ 
{\frac {1-x^{(d+1)l}}  {1-x^{dl + 1}}}\right )^{1/2}\ 
 e_{J_l} s_{J,l-1} t.
$$
\end{lemma}

\v\ni
\begin{proof}
It follows by induction from Lemma \ref{Lemma 3.2.2} that
$$
e_{J_l} e_{J, l-1} t = 
\prod_{i = 1}^{l-1} (a_{\pm(dl + l-i)})^{i/2}\  e_{J_l} s_{J,l-1} t,
$$
so the result follows  from an application of  Lemma \ref{Lemma 3.2.3}. 
\end{proof}

\v
Write $P^t$ for the coefficient in Lemma \ref{Lemma 3.2.4}.
Now suppose that $t$ and $t'$ are two $k$-special $l$-straight
skew tableaux with the
same initial shape satisfying conditions (1) and (2) above.  Then
$$\langle  t', Et\rangle = 
\langle  e_{J, l-1} t', e_{J_l} e_{J, l-1} t\rangle = 
P^t P^{t'} \langle s_{J, l-1}t', e_{J_l} s_{J, l-1} t \rangle.
$$
The skew tableau $ s_{J, l-1}t$ has the elements of 
$J_l$ in increasing order in row $a$ and those of 
the complementary set $\{r + 1 : r \in J_l\}$
in increasing order in row $b$.  Similarly  $ s_{J, l-1}t'$ has
the elements of $J_l$ in some row $a'$ and those of the complementary set
in some row $b'$.  Evidently  
$\langle s_{J, l-1}t', e_{J_l} s_{J, l-1} t \rangle = 0$ unless
$\{a, b\} = \{a', b'\}$, and either $t = t'$, or $ t' = s_E t$.
One has 

\v\ni
\begin{proposition}
\label{Proposition 3.2.5}
$$\langle E t,  t\rangle =( P^t)^2 a_{dl}^l = 
x^{\binom l  2}\  {\frac{1- x^{(d+1)l}} {1-x^{dl}}} = x^{\binom l 2} a_d(x^l),
$$
and
$$\langle   Et, s_E t\rangle = P^t P^{s_E t} a_{\pm ml}^{l/2} = 
x^{\binom l  2}\  (a_d(x^l) a_{-d}(x^l))^{1/2}.  \qed
$$
\end{proposition}

Note that the matrix coefficients of $x^{- \binom l 2} E$ with respect to
$k$-special $l$-straight tableaux are just matrix coefficients of 
$e_1(x^l)$ in the orthogonal representation. (Compare section 1.3).
\v
\ni
\subsection{An embedding of the $k$-row quotient of $H_m(x^l)$.}
\v
Given an natural number $m$, we consider a sequence of big diamond
elements $E_1, E_2, \dots, E_{m-1}$ in a certain $H_n = H_n(x, K(x))$, as follows. 
Let
$n = n_0 + ml$, where $n_0 = (l-1) {\binom{k}{2}}$
denotes the size of the smallest
$k$-critical Young diagram $(l-1)\rho$.

Let $E_i$ denote the big diamond element in $H_n$,
$$E_i = E(n_0 + (i-1)l),$$ for $1 \le i \le m-1$.  The corresponding
permutation $s_{E_i}$ affects, at most, the integers between
$n_0 + (i-1)l$ and $n_0 + (i+1)l$. 

Consider $H_n$ acting in the orthogonal
representation on tableaux of length $n$; for any  tableau $t$,
$E_i t$ is a linear combination of tableaux $s$ with $s(j)=t(j)$
for $j\leq n_0 + (i-1)l$ or $j\geq n_0 + (i+1)l$;  in particular,
$s(j)=t(j)$ for $j \le n_0$.   

Fix a tableau $t$ of length $n$ satisfying the following conditions:
\hfill\break
\vbox{
\smallskip
(1)\ 
$t$ passes through $(l-1)\rho$.

\smallskip (2)\ 
$t$ is a $k$-special  $l$-straight skew tableau between $(l-1)\rho$ and
its final diagram $\la(t)$.
}

\smallskip
One can check that for any tableau $\tilde t$ which
is $l$-equivalent to $t$ and whose final shape $\la(\tilde t)$ has
at most $k$ rows, $\tilde t$ also
satisfies conditions (1) and (2). Conversely, any tableau satisfying
conditions (1) and (2) is $l$-equivalent to $t$.

We recall several notations: 
The sum of all path idempotents $p_s$, where
$s$ is equivalent to $t$ is denoted by $p_{[t]}$,
see Definition \ref{Definition 1.5.2} and
Lemma \ref{Lemma 1.5.2}.
The direct sum of the $H_n(x, K(x))$ modules $V_\la$, where $\la$ ranges over
Young diagrams of size $n$ with at most $k$ rows is denoted by
$V(n,k)$; the representation of $H_n(x, K(x))$ on $V(n,k)$ is denoted by $\Pi_k$.
The representation of $H_n(x, K(x))$ on $(K(x)^k)^{\otimes n}$ with respect to
the natural basis is denoted by $\Phi_k$.

Note that the range of 
$\Pi_k(p_{[t]})$ 
is spanned by the set of  standard tableaux satisfying conditions (1)-(2).

\v\ni
\begin{theorem}
\label{Theorem 3.3.1}
\medskip\ni
\begin{enumerate}
\item
The images of $p_{[t]}$, $p_{[t]}E_1 p_{[t]}$, \dots,
$p_{[t]}E_{m-1} p_{[t]}$ in the $k$-row quotient of $H_n(x, K(x))$
generate a subalgebra isomorphic to the $k$-row quotient of
$H_m(x^l, K(x))$.
\item
The range of $\Pi_k(p_{[t]})$ in  
$V(n,k)$
is isomorphic to 
${\mathcal R}  \otimes V(m,k)$,
where $\mathcal R$ is the $K(x)$-vector space with basis
 the set of tableaux which are
$l$-equivalent to $t_{|(l-1)\rho}$.  Furthermore, 
 $x^{-{\binom l  2}}p_{[t]}E_i p_{[t]}$ acts as $\ident \otimes e_i(x^l)$ on 
${\mathcal R}  \otimes V(m,k)$.
\item
Let $q \in K$ be a  primitive $l$-th root of $1$.   The images of 
$p_{[t]}(q)$, $p_{[t]}E_1 p_{[t]}(q)$, \dots,
$p_{[t]}E_{m-1} p_{[t]}(q)$ in the $k$-row quotient of 
$H_n(q)$ generate a subalgebra isomorphic to the $k$-row quotient
of $KS_m$.
\end{enumerate}
\end{theorem}
\v
\ni
\begin{proof} Consider $H_n(x, K(x))$ acting by the orthogonal
representation
$\Pi_k$
on
$V(n,k)$.
 For a tableau $s$
$l$-equivalent to $t$,  consider
 $s_{|(l-1)\rho} \otimes \Psi(s)$, where
$\Psi$ is as in Lemma \ref{Lemma 3.1.1}.  The map
$$\Theta : s \mapsto
s_{|(l-1)\rho} \otimes \Psi(s)$$
 extends to a linear isomorphism
from the range of $\Pi_k(p_{[t]})$ onto ${\mathcal R} \otimes
V(m,k)$, according to the
remarks preceding the statement of the Theorem.

The formulas of Proposition \ref{Proposition 3.2.5}
show that
$$\Theta \circ x^{-{\binom l  2}}  \Pi_k(p_{[t]}E_i p_{[t]}) \circ \Theta\inv
= \ident \otimes \Pi_k(e_i(x^l))$$
acting on ${\mathcal R} \otimes V(m,k)$.
This proves points (a) and (b), since
$\Pi_k(H_n(x, K(x)))$ is isomorphic to the $k$-row quotient of 
$H_n(x, K(x))$.

For part (c),  consider the elements $\tilde T_i = x^{-{\binom l  2}}
 p_{[t]}E_i p_{[t]}  + x^l p_{[t]}$ in $H_n(x, K(x))$, which are evaluable
at the $l$-th root of unity $q$.

It follows from part (b),  and the isomorphism
$\Pi_k(H_n(x, K(x))) \cong \Phi_k(H_n(x, K(x)))$, 
 that the assignment $T_i \mapsto \Phi_k(\tilde
T_i)$ determines a non-unital $K(x^l)$-algebra homomorphism
$$
\varphi : H_m(x^l, K(x^l)) \rightarrow \Phi_k(H_n(x, K(x))),
$$
with the  image of $\varphi$ 
isomorphic to the $k$-row quotient of $H_m(x^l, K(x^l))$.
Since $\Phi_k(\tilde T_i)$ is evaluable at $x= q$, the representation
$\varphi$ is evaluable at $x = q$\   (that is, at $x^l = 1$).
By Lemma \ref{Lemma 0.4.1}, $\varphi$ induces a $K$-algebra homomorphism
$$
\tilde\varphi :KS_m \rightarrow \Phi_k(H_n(q, K)),
$$
satisfying $\tilde \varphi(s_i) = \Phi_k(\tilde T_i)(q)$.
We have to show that the image of $\tilde\varphi$ is isomorphic to the
$k$-row quotient of $KS_m$.

Let $s$ be a standard tableau of size $m$, and let $p_s$ be
the corresponding path idempotent in $H_m(x^l, K(x^l))$.
According
to Lemma \ref{Lemma 1.5.4}, $p_s$ is evaluable at $x^l = 1$, and its
evaluation
${p_s}_{|\, x^l = 1}$ is the corresponding path idempotent in $KS_m$.
We need to show that $\tilde\varphi({p_s}_{|\, x^l = 1})$ is non-zero if, and
only if, $s$ has no more than $k$ rows.

Observe that $\varphi(p_s)$ is non-zero if, and only if $s$ has no more
than $k$ rows, and  as
$\varphi(p_s)$ is an idempotent, its rank is the same as that of
$\varphi(p_s)_{|\, x= q}$.  But
$\varphi(p_s)_{|\, x= q} = \tilde\varphi({p_s}_{|\, x^l = 1}) $
by Lemma \ref{Lemma 0.4.1}, so    the desired conclusion follows.
\end{proof}

\v\ni
\begin{lemma} \label{Lemma 3.3.2} Let $y$ be a point in $D^+\cap \Z^k$, and
 let $0\leq i,m\leq k$ be integers such that $i+m\leq k$. 
Assume that $y$ has $i+m$ coordinates which
are congruent to $ y_1 \ \mod  l$,
among them the first $i$ (i.e. $y_1$ up to $y_i$). Let 
$e^{(i)}=(1,\ ...,\ 1,0,\ ..., 0)$ (with $i$ 1's and $k-i$
0's), and let, for given $i\in \N$, $t_\rho$ be the path from $y$ 
to $y+e^{(i)}$ such that $t_\rho(j)=y+e^{(j)}$, for $j=1,2,\ ...\ i$.
\begin{enumerate}
\item
Any path $\tilde t_\rho$ which is equivalent to $t_\rho$ 
ends in $y+e^{(i)}$ or a point $\leq y+e^{(i)}$ in lexicographic order.
\item
There exist at most $(m+i)!/m!$ paths in $D^+$ which are equivalent to 
$t_\rho$, among which $i!$ paths end in $y+e^{(i)}$.
\item
Let $w\in W$ such that $w(t_\rho)$ is a path corresponding
to a skew tableau. Then there exist at most
$(m+i)!/m!$ paths starting from $w(y)$ which are equivalent to
$w(t_\rho)$.
\end{enumerate}
\end{lemma}
\v
\ni
\begin{proof}
Assume, for simplicity, that $y_1\equiv 0 \ \mod  l$. Then
$t_\rho(j)\equiv 1 \ \mod  l$, for $j=1,2,\ ...\ i$, and the same
has to hold for any tableau $\tilde t_\rho$ which is equivalent to $t_\rho$.
For such a $\tilde t_\rho$, we therefore have $m+i$ possibilities
for the first box, $(m+i-1)$ possibilities for the second one, and so on.
Altogether, there are $(m+1)(m+2)\ ...\ (m+i)$ possibilities, as stated.
There are $i!$ different ways of adding a box in each of the first $i$
rows. This shows (b). The proof of (c) goes similarly. (a) is clear. 
\end{proof}

\v\ni
\begin{lemma} \label{Lemma 3.3.3} 
Let $K$ be a field of characteristic 0, and let 
$q$ be a primitive $l$-th root of unity in $K$. 
Let $\la$ be a  $k$-critical diagram. 
Then there exists
an evaluable idempotent $p\in H_n(x, K(x))$ such that $pS^\la \ne 0$, and
for all diagrams $\mu$ such that    $\mu\neq \la$ and 
$\mu$ has at most $k$ rows,
$pS^\mu = 0$.
\end{lemma}

\v
\ni
\begin{proof}
Let us assume first that $\la$ 
is a reduced $k$-critical diagram, 
i.e. that $\la_k=0$. Let $t$ be a tableau of shape $\la$ 
satisfying conditions (1) and (2) above.
Using the formulas for central idempotents
of $K S_m$ and Theorem \ref{Theorem 3.3.1}, we obtain
an evaluable element $\tilde p$
in $p_{[t]}H_np_{[t]}$ such that $\pi_\mu(\tilde p)$ is evaluable
at $q$ and an idempotent for all $k$-row diagrams $\mu$, and  is nonzero
if and only if $\mu = \la$. By Proposition \ref{Proposition 0.7}, there exists
an evaluable idempotent $p\in H_n$ such that for diagrams $\mu$
of at most $k$ rows, $p$ is  nonzero
on $S^\mu$  if and only if $\mu = \la$.

The general case is proved by induction on $\la_k$, the last component
of $\la$. The case $\la_k = 0$ has already been verified.
For $\la_k>0$, we can assume that the
claim has been verified for the diagram $\la-e^{(k)}$;
more specifically, we can make the inductive hypothesis that
 there exists
a tableau $t$ of shape $\la-e^{(k)}$ and 
an evaluable subidempotent $p$ of $p_{[t]}$ which is nonzero only
on $S^{\la}$. It follows from Lemma \ref{Lemma 3.3.2} that for any 
extension $\hat t$ of $t$ into $\la$, all its equivalent tableaux
going through $\la-e^{(k)}$ also end up in $\la$ (as $r=m=k$).
Hence $pp_{[\hat t]}$ has all the required properties. 
\end{proof}
\v
The following result is contained in results of James and Dipper
([DJ, 4.11) and James and Mathas [JM].
\v\ni
\begin{corollary} \label{Corollary 3.3.4}  
Let
$q$ be a primitive $l$-th root of unity in $K$. 
For any $k$ and
for any $k$-critical diagram $\la$, $S^\la = D^\la$.  Furthermore,
$d_{\nu \la} = 0$ for all $\nu \ne \la$ such that $\nu$ has
at most $k$ rows.
\end{corollary}
\v\ni
\begin{proof} Suppose one has established that $S^\la = D^\la$ for a 
{\em particular} $k$-critical diagram $\la$.
By Lemma \ref{Lemma 3.3.3}, there is an 
evaluable idempotent $p \in H_n(x)$ with the property that 
$p$ is non-zero on $S^\la$ but zero on any
$S^\mu$, if  $\mu\neq \la$ and
$\mu$ has at most $k$ rows.
It follows from Proposition \ref{Proposition 1.7.1} that $d_{\nu \la} = 0$ for
all diagrams $\nu \ne \la$ with at most $k$-rows.

Now we proceed by induction on the reverse
lexicographic order on $k$-critical diagrams  of fixed size $n$.
Consider  the diagram $\la$
which is largest in lexicographic  order among $k$ critical diagrams
of size $n$.
Then $\la$ is clearly the highest diagram in its $W$ orbit,
and hence  $D^\la=S^\la$, by Lemma \ref{Lemma 2.6}.
Suppose that the claim has been verified for all $k$-critical
diagrams $\tilde \la$ such that $|\tilde \la| = |\la| = n$ and
$\tilde \la > \la$.  By the Nakayama conjecture, if 
$\mu$ is a diagram such that $D^\mu$ is a composition factor
of $S^\la$, then $\mu \trianglerighteq \la$, so in particular
$\mu$ has at most $k$ rows, and $\mu$ is a $k$-critical diagram.
But then it follows from the induction hypothesis that if $\mu \ne \la$,
then $D^\mu = S^\mu$, and $d_{\la \mu} = 0$.  Therefore
$S^\la$ has no composition factors other than $D^\la$, so
$S^\la = D^\la$. 
\end{proof}
\v\ni
\begin{remark} \label{Remark 3.3.5} If one knows {\em a priori} that 
$S^\la = D^\la$  for {\em all} $k$-critical diagrams, then
it follows at once from the Nakayama conjecture that 
$d_{\nu \la} = 0$ for all $k$-critical diagrams $\la$ and for
all diagrams $\nu \ne \la$ with at most $k$ rows.
\end{remark}
\v\ni
\subsection{General $l$-straight tableaux}
There is a version of 
Proposition \ref{Proposition 3.2.5} which holds for arbitrary (i.e. not $k$-special)
$l$-straight skew tableaux.  The matrix coefficients are  computed
only modulo $q^l = 1$.  The idea of the computation is similar to that 
of Proposition \ref{Proposition 3.2.5}, but the combinatorics are more complicated. 
We  omit the details and merely state the result.  We use this result
in Section 5 for our remarks on the ``boundary region".

\v\ni
\begin{proposition}\label{Proposition 3.4} Suppose $t$ and $t'$ are $l$-straight 
skew tableaux of length $2l$ with the
same initial and final shapes. Let $Q^{t,t'}=\langle t',Et\rangle$,
and let $Q^t=Q^{t,t}$. Then 
\smallskip\ni
\begin{enumerate}
\item
$Q^t(q) =$ $\displaystyle C(q)
{\frac{(d(t;l,l+1)+1)\prod_{i=1}^{l-1} (d(t;i+1,l+i)-1) }
 {\prod_{i=1}^l d(t;i,l+i)}}$, where $C(q)$ does not depend on $t$.
\item
If
there is a subset $A \subseteq \{1,2,\dots,l\}$,
such that $\displaystyle t' = 
\left [ \prod_{a\in A}(a,a+l)\right ] t$, then
$Q^{t,t'}(q) = \sqrt{Q^t(q)Q^{t'}(q)}$;
otherwise $Q^{t,t'}(q)=0$.  Here $(a, a+l)$ denotes the transposition
which interchanges $a$ and $a+l$.
\end{enumerate}
\end{proposition}
\v\ni
\begin{remark}The condition in (b) is equivalent to the existence of
a subset $L \subseteq J_l$ such that  $s_J t' = s_L s_J t$.
\end{remark}

\section{Estimates on the decomposition numbers $d_{\la \mu}$.}

In this section we obtain bounds on the decomposition numbers
$d_{\la \mu}$ for the Hecke algebra $H_n(K,q)$ where $K$ is a field
of characteristic 0, and $q$ is a primitive $l$-th root of unity.

\v\ni
\subsection{Critical points associated to interior points.} We call a point
$y\in D^+$  {\em interior} if $y_i-y_{i+1}\geq l$ for
$i=1,2,\ ...\ k-1$. We want to access an interior point $y$
from a convenient $k$-critical point.

 Write $\Mod(n,l) = n - [n/l]l$; thus $n = [n/l]l + \Mod(n,l)$, and
$0 \le\Mod(n,l) \le l-1$
\v\ni
\begin{definition}
Let $y$ be an interior point with integer
coefficients. Its {\em associated
critical point} $c=c(y)$ is determined as follows:
Put $r_k = 0$ and $r_i = \Mod(y_i - y_{i+1}, l)$ for $1 
\le i \le k-1$.  Set $c_i = y_i - \sum_{j \ge i} r_j$ for $1 
\le i \le k$. 
\end{definition}
\v\ni
\begin{lemma} \label{Lemma 4.1.1}
Let $y$ be an interior point, and let
$c$ be its associated critical point. 
\v\ni
\begin{enumerate}
\item
$c$ is a $k$-critical interior point.
\item
$y-c = \sum_{i=1}^{k-1}r_i e^{(i)}$, and $|y| -|c| = \sum_{i=1}^{k-1}i\ r_i
\le (l-1)k(k-1)/2$.
\item
There exists a path
$t_\rho$ from $c$ to $y$  such that all its conjugates
end up in $y$ or in lexicographically lower diagrams. The number of its conjugates
is less or equal to
$$\prod_{i=1}^{k-1}\prod_{j=d_{i+1}}^{d_i-1} {\frac{(m(j,i)+i)!} {m(j,i)!}},$$
where $m(i,j)=|\{ y_s : \ y_s\equiv j \ \mod  l,\ s>i\}|$. 
Among those conjugate paths, exactly $\prod_{i=1}^{k-1} (i!)^{r_i}$
end in $y$.
\end{enumerate}
\end{lemma}
\v\ni
\begin{proof}
One has $c_i - c_{i+1} = y_i - y_{i+1} - r_i =
[(y_i - y_{i+1})/l]\ l \ge l$, since $y_i - y_{i+1} \ge l$.  This shows
(a).  Point (b) is evident from the definitions; the last inequality
results from $r_i \le l-1$.

For (c) we take the path $t_\rho$ from $c$ to $y$,
\begin{equation}
\begin{split}
y^{(k)} = c &\rightarrow y^{(k-1)} = y^{(k)} + 
r_{k-1} e^{(k-1)} \rightarrow \dots \cr
&\rightarrow		y^{(i)} = y^{(i+1)} + r_{i} e^{(i)}			\rightarrow \dots \cr
&\rightarrow		y^{(1)} = y,
\end{split}
\end{equation}
where one goes from $y^{(i+1)}$ to $y^{(i)}$ by adding the $(i\ r_i)$ cells
column-wise, as in Lemma \ref{Lemma 3.3.2}.  (Note that $y^{(i+1)} +s e^{(i)}$ has
at least the first $i$ rows congruent ($\mod\  l$), for
$0 \le s \le r_i$.	

If $\tilde t_\rho$ is another path of length $\sum i\ r_i$ starting from
$c$ which is equivalent to $t_\rho$, then by the definition of 
equivalence, the two paths are equivalent at each step.  Therefore, it
follows from Lemma \ref{Lemma 3.3.2} and induction that $\tilde t_\rho$ ends in $y$ or in
a point which is lexicographically lower than $y$.  
The estimates about the number of all equivalent paths,
and the number of those which end in $y$ follows by induction
on $|y| - |c|$, using Lemma \ref{Lemma 3.3.2},(b) and (c). 
\end{proof}

\v
\subsection{An estimate on $\dlm$}
Let $c$ be an interior $k$-critical point and $t_\rho$ a path from
$c$ to some interior point $y$.  For any diagram $\la$ with
at most $k$ rows, denote by 
$N(t_\rho, \la)$ the number of paths from $c$ to $\la + \rho$ which are
equivalent to $t_\rho$.

\v\ni
\begin{theorem}
\label{Theorem 4.2}
Let $\mu$ be a diagram such that 
$y = \mu + \rho$ is an interior point.
Then $D^\mu\neq 0$. 
Let $c$ be the critical point associated to
$\mu$ and $t_\rho$ the path from $c$ to $\mu + \rho$ constructed in
Lemma \ref{Lemma 4.1.1}.   If $\la$ is a
diagram with $\leq k$ rows, the multiplicity $d_{\la\mu}$ of $D^\mu$ in $S^\la$
satisfies
$$d_{\la\mu} \le N(t_\rho, \la)/N(t_\rho, \mu).$$
\end{theorem}
\v\ni
\begin{proof}
Since $\mu + \rho$ is an interior point, it follows that $\mu$ has
no two rows of the same length; so $\mu$ is $l$-regular and $D^\mu \ne 0$.

 By  Corollary \ref{Corollary 3.3.4},
there exists an evaluable idempotent $z_{(c-\rho)}$ which acts
as the identity on $S^{c-\rho}$ and as zero on any Specht module $S^\nu$
belonging to a Young diagram $\nu$ with $k$ rows at the
most such that 
$|\nu| = |c - \rho|$, and $\nu \ne c-\rho$. Let $t^{(0)}$ be any tableau of shape
$c-\rho$, and  let $\hat t_\rho$ be an extension of $t^{(0)}_\rho$ into $y$ 
by the  path $t_\rho$  constructed in
Lemma \ref{Lemma 4.1.1}. Then $z_{(c-\rho)}p_{[\hat t]}$ is an evaluable
idempotent in $H_n$ ($n = |\mu|$).  

$\Phi_k(z_{(c-\rho)}p_{[\hat t]})$ is equivalent to the projection 
(in $\bigoplus_{\ell(\la) \le k} V_\la$)
on  the span of all those  paths which are equivalent to
$\hat t_\rho$ and which go through $c$.
Hence its evaluation acts
as zero on any $S^\la$ with $\la > y-\rho=\mu$, but not as zero on $S^{\mu}$.
It follows that the rank $r$ 
of $p=z_{c-\rho}p_{[\hat t]}$ on $S^\mu$ is also the rank by which it acts on
$D^\mu$. 

The rank $r$ is equal to the number of paths which end in $y$ and go through $c$.
If $r_0$ is the number of paths which are equivalent to $t_\rho^{(0)}$
and end in $c$, one sees that $r=r_0N(t_\rho ,\mu)$.
Similarly, the rank by which $p$ acts on $S^\la$ is equal to
$r_0N(t_\rho ,\la)$. By Proposition \ref{Proposition 1.7.1}(b) 
\begin{equation*}
d_{\la\mu}
\leq (r_0N(t_\rho ,\la))/(r_0N(t_\rho ,\mu)) = 
N(t_\rho ,\la))/N(t_\rho ,\mu). \quad 
\end{equation*} 
\end{proof}

\v\ni
\subsection{Reduced paths.}  In the bound for $d_{\la \mu}$ in
4.2, the denominator is known, namely 
$$N(t_\rho, \mu) = 
\prod_{i = 1}^{k-1} (i !)^{r_i}.$$  
The numerator is divisible by
this quantity, and it remains to describe the quotient.
 
We define the {\em reduced path} $t_{\rho ,{\rm red}}$ of $t_\rho$ to be the
sequence of points:
\begin{equation}
\begin{split}
c, c + e^{(k-1)}, c + 2e^{(k-1)},&\dots, c + r_{k-1}e^{(k-1)} = y^{(k-1)}, \cr
y^{(k-1)}, y^{(k-1)} + e^{(k-2)}, &\dots, y^{(k-1)} + r_{k-2} e^{(k-2)} =
y^{(k-2)}, \cr 
\quad\quad&\dots \cr
y^{(2)},y^{(2)} + e^{(1)}, &\dots, y^{(2)} + r_1 e^{(1)} = y.
\end{split}
\end{equation}
That is, we divide the path $t_\rho$ into segments, the 
first $r_{k-1}$ of length $k-1$, the next $r_{k-2}$ of length $k-2$,
and so forth, until the last $r_1$ segments of length $1$; the reduced
path is the list of endpoints of these segments.
Similarly, if $\tilde t_\rho$ is a path equivalent to $t_\rho$, we divide
this path into segments of the same lengths, and define the 
reduced path $\tilde t_{\rho,{\rm red}}$ to be the sequence of endpoints
of these segments.

We define $n(\la, \mu)$ to be the number of distinct
reduced paths belonging to paths $\tilde t_\rho$ which are equivalent to
$t_\rho$ and which end in $\la + \rho$.

\v\ni
\begin{lemma} \label{Lemma 4.3.1}
Let $\mu$ be a diagram such that 
$y = \mu + \rho$ is an interior point.
Let $c$ be the critical point associated to
$y$ and $t_\rho$ the path from $c$ to $y$ constructed in
Lemma \ref{Lemma 4.1.1}. 
Let $\la$ be a diagram with $k$ rows at the most. Then 
$$
N(t_\rho ,\la)/N(t_\rho ,\mu) = n(\la, \mu)
$$ 
\end{lemma} 
\v
\ni
\begin{proof}It follows from Theorem \ref{Theorem 4.2} that there is only one reduced
path which ends in $y$.  It suffices to observe that if $\la$ is
a diagram such that $N(t_\rho ,\la) \ne 0$, then the ratio
$N(t_\rho ,\la)/n(\la, \mu)$ is the same as
$N(t_\rho ,\mu)/n(\mu, \mu) = N(t_\rho ,\mu)$. 
\end{proof}

\v\ni
\begin{corollary} \label{Corollary 4.3.2}
Let $\mu$ be a diagram such that 
$y = \mu + \rho$ is an interior point.
Let $c$ be the critical point associated to
$y$ and $t_\rho$ the path from $c$ to $y$ constructed in
Lemma \ref{Lemma 4.1.1}.   If $\la$ is a
diagram with $\leq k$ rows, the multiplicity $d_{\la\mu}$ of $D^\mu$ in $S^\la$
satisfies
$$
d_{\la \mu} \le n(\la, \mu)
$$
\end{corollary}
\v\ni
\begin{remark} \label{Remark 4.3.3}
The reduced paths described above can also be
interpreted as certain semi-standard tableaux of shape $\la  \setminus
(c- \rho)$.   Namely, let $\nu$ be the conjugate diagram
of $(\mu + \rho -c )$, that is, \break
 $\nu = (k-1)^{r_{k-1}}(k-2)^{r_{k-2}}\cdots 1^{r_1}$.
Then the reduced paths correspond 1-to-1 with semi-standard tableaux
of of shape $\la  \setminus
(c-\rho)$ and weight $\nu$  (that is, with $\nu_1$ 1's, $\nu_2$ 2's, etc.)
such that each entry of $j$ has content  (column index  minus
row index) congruent to $j-k \ \mod  l$. 
\end{remark} 
\v\ni
\begin{remark} \label{Remark 4.3.4}
 Observe that the hyperplanes passing through $c$ generate
a reflection group $W_c$ which is isomorphic to $S_k$. 
Consider the set $\A$ of all alcoves through which $t_\rho$ and its
conjugates run. Let $A\in \A$. Assume that 
for all
$w\in W_c$, $w(A)$ is an alcove
in the Weyl chamber which does not touch the boundary of the Weyl chamber. 
Then $n(\la ,\mu)=n(w(\la +\rho)-\rho ,\mu)$ for all
$w\in W_c$. We conjecture that a similar symmetry also holds
for the $d_{\la\mu}$.
\end{remark}
\v\ni

\begin{corollary}
\label{Corollary 4.4}
$$
\sum_{\la ,\ \la_{k+1}=0}\ d_{\la\mu}\leq
\prod_{i=1}^{k-1}\ \prod_{j=d_{i+1}}^{d_i-1}
\begin{matrix}
m(j,i)+i \cr 
i\cr
\end{matrix},
$$
where $m(i,j)=|\{ y_s,\ y_s\equiv j \ \mod  l,\ s>i\}|$. 
\end{corollary}
\v
\ni
\begin{proof}
The right hand side is less or equal to the number of all paths
within $k$ row diagrams which are equivalent to $t_\rho$ divided by the number
of paths equivalent to $t_\rho$ which end  in $y=\mu+\rho$. 
Our estimate follows from Lemma  \ref{Lemma 4.1.1}(c). 
\end{proof}

\v
\ni
\section{Examples}

We give a geometric description of our results for $k=3$.
The  weights $L_i$ ($1 \le i \le 3$) of ${sl}_3$ in the vector
representation ($L_3 = -L_1 - L_2$) can be represented by  
three coplanar unit vectors which are mutually equidistant.

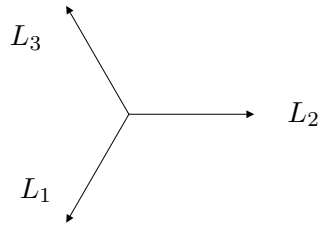
\centerline{
{\rlap{\raise .4 in \hbox{$L_3$}}}
{{\raise -.4 in \hbox{$L_1$}}}
\BoxedEPSF{weights.eps scaled 300}
{{\kern .1 in \hbox{$L_2$}}}
}

Consider the map $\R^3 \rightarrow \R^2$, 
$ y \mapsto \sum y_i L_i$.
Partitions $\la$ of length $\le 3$  map to 
dominant integral weights $\la_1 L_1 + \la_2 L_2 + \la_3 L_3$.
The open Weyl chamber is the image of $D^+$, namely
$C = \{\sum y_i L_i : y_1 > y_2 > y_3 \ge 0\}$.  A path $t_\rho$
corresponding to a tableau or skew tableau $t$ is mapped to
a path in $C$,
each of whose segments is one of the vectors $L_i$ and whose endpoints
are integer points of $C$, namely in the integer span of the
$L_i$.  $C$ is divided into triangular tiles by the lines
$\{\sum y_i L_i : y_i - y_j = m l\}$.  The $3$-critical points of $D^+$ map
to  the intersection of three such lines.  The boundaries of $C$ are
formed by the lines $\{\sum y_i L_i : y_2 - y_3 = 0\}$ and
$\{\sum y_i L_i : y_1 - y_2 = 0\}$.  For any point $y$, we define the
closed positive cone of $y$ to be $y + \bar C$.
The interior points of $D^+ \cap \Z^3$ map to the cone
$2l L_1 + l L_2 + \bar C$.

\smallskip
\centerline{\BoxedEPSF{weyl.eps scaled 800}}

\smallskip
In the following statements, we no longer distinguish between points
of $D^+$ and the corresponding points of $C$.
Let $y$ be an interior point which does not lie on any of the
affine hyperplanes (i.e. lines in our case) just described.
Let $c$ be its critical point.
The following is easy to check:

(a) $(y-c)_1<2l$ and $(y-c)_1\neq l$,

(b) $(y-c)_1<l$ if and only if $y$ is in a triangle pointing
upwards, with its highest vertex equal to $c$ (see next picture
below),

(c) $(y-c)_1>l$ if and only if $y$ is in a triangle pointing
downwards; here $c$ is the highest vertex of the triangle
above the one containing $y$ (see second picture below).

The following two diagrams show  a path from $c$ to $y$ together with
its conjugate paths in the two cases:

\smallskip
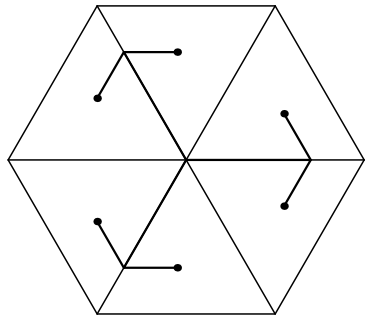
\centerline{\BoxedEPSF{hex.eps scaled 500}}
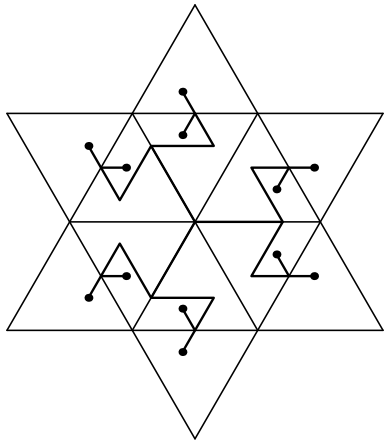
\smallskip

\centerline{\BoxedEPSF{david.eps scaled 600}}

\smallskip
Next we consider the boundary region (that is the complement of
the interior region) in $D^+$ for $k=3$.  This can be
extended to certain boundary points also for larger $k$
(see Remark 3 below). Let $y$ be a point
in the left boundary region, i.e. such that $y_2 - y_3 < l$.
We can assume $y_3=0$, and hence $y_2<l$.
We can therefore approach $y$ by a path as shown in the figures
below, which also show the conjugate paths.

\smallskip
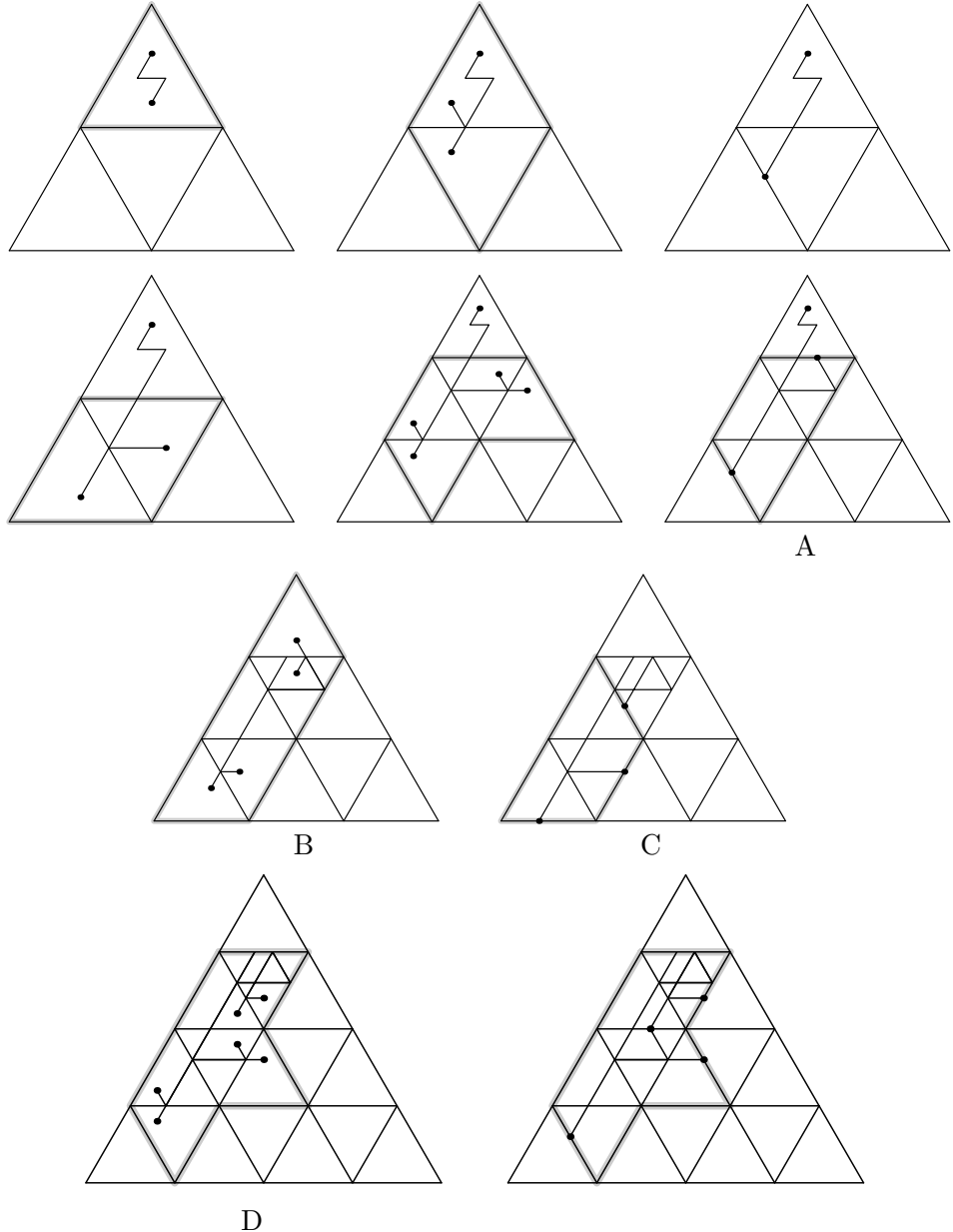
\centerline{\tBoxedEPSF{step1.eps scaled 400 }\quad\tBoxedEPSF{step2.eps
scaled 400 }\quad \tBoxedEPSF{step3.eps scaled 400 }}

\smallskip
\centerline{\tBoxedEPSF{step4.eps scaled 400 }\quad\tBoxedEPSF{step5.eps
scaled 400 }\quad 
\tBoxedEPSF{step6.eps scaled 400 }\rlap{ \kern -.9 in \raise -1.5 in
\hbox{A}}}

\smallskip
\centerline{\tBoxedEPSF{step7.eps scaled 400 }
\rlap{ \kern -.9 in \raise -1.5 in
\hbox{B}}
\quad\tBoxedEPSF{step8.eps
scaled 400 }
\rlap{ \kern -.9 in \raise -1.5 in
\hbox{C}}
}

\smallskip
\centerline{\tBoxedEPSF{step9.eps scaled 500 }
\rlap{ \kern -1.2 in \raise -1.9 in
\hbox{D}}
\quad\tBoxedEPSF{step10.eps
scaled 500 }
}

\smallskip
Here comes the crucial step for such paths:
Consider the picture to the right of picture D.
Our path orbit contains 5 paths. The corresponding
projection $p_{[t]}$ can be split into the sum of 2
subidempotents as follows:
It is easy to check that the endpoint $a$ of the lowest and
the endpoint $b$ of
the rightmost path are the alphabetically highest and
secondhighest points among all the endpoints of paths in that picture.
Applying Proposition \ref{Proposition 3.4} to the skew tableaux containing the last
$2l$ steps of our paths, we obtain an evaluable
 subidempotent $p$ of $p_{[t]}$
which acts as 0 on the lowest path and as 1 on the rightmost
path. It again follows from Proposition \ref{Proposition 3.4} that $p$ acts
as a rank 3 idempotent on the space spanned by our 5 paths.
Hence $p_{[t]}-p$ can be described by picture $E$.

\smallskip
\centerline{\tBoxedEPSF{step11.eps scaled 500}
\rlap{ \kern -1.2 in \raise -1.9 in
\hbox{E}}
}

\v
This is not to be taken quite literally; namely $p_{[t]}-p$ is not the
sum of two ordinary path idempotents.  Nevertheless, we are essentially
back to the pattern of figure A in the series, and as
the path $t_\rho$ is extended into further alcoves the sequence
of figures A through E repeats itself.

\v\ni
{\bf Remarks}
1. It is possible to prove a similar pattern for the right boundary region.
Essentially, we use the mirror path of the path constructed for the
left boundary region. Such a path can be interpreted as obtained
by tensoring by the dual representation of the fundamental representation
of $sl_3$; in path language, it would mean we construct all path
extensions of length 2, and let the projection $p_{[1^2]}$ act
on these extensions. It can be shown that the same algebraic
theory can be applied as for the ones studied before, e.g.
by using the result in [KW] which characterizes fusion categories
of type $A$.

2. For $k=3$ all our bounds are sharp, as can be seen by comparing the results with
computations using the algorithm of [LLT].  However, for $k \ge 4$, our bounds are
no longer sharp.

3. Our method for dealing with boundary points for $k=3$
can be generalized to get partial results for $k \ge 4$, but it is more
efficient to use the LLT algorithm.

\v
\ni
\section{Connection to quantum groups}
\subsection{} Let $\U$ be the Drinfeld-Jimbo quantum group corresponding
to the root system $A_{k-1}$. In the following we will
also assume $q$ to be a primitive $l$-th
root of unity with $l>k$; in this case we take for $\U$
the version of the quantum group
as defined by Lusztig, where $q=v^2$ (see [Lu]). 

The notion of {\em tilting modules}
of a quantum group was introduced by Andersen [A], inspired by a similar
notion for algebraic groups which was defined by Donkin [Do1].
For type $A$, it can be shown that any direct summand of $V^{\otimes n}$
is a tilting module,
where $V$ is the fundamental $k$-dimensional $\U$ module. This follows
from the fact that
$V$ is a tilting module, and tensor products and direct summands
of tilting modules are also tilting; see [A] for details.
Conversely, it is possible to characterize tilting modules for type A
as  direct sums of direct summands
of $V^{\otimes n}$s (with $n$ varying).

It was already observed by Jimbo [Ji] that the image of the Hecke algebra
under the representation on $V^{\otimes n}$, as defined in Section 2,
commutes with the action of $\U$, and, if $q$ is not a  proper root of
unity,   the image of the Hecke algebra and the image of the quantum group
are commutants of one another; this is a $q$-analog of the famous
Schur-Weyl duality between the general linear group and
the symmetric group. 

It has only recently been established that the same duality holds,
between the Lusztig version of the quantum group of type $A$ and the Hecke 
algebra, 
when $q$ is a root of unity [DPS].
(See [J] for the case of algebraic groups
with positive characteristic). It follows that indecomposable tilting modules 
for $\U$ all have the form $p V^{\otimes n}$, for some $n$ and for some
primitive idempotent $p$ in the $k$-row quotient of $H_n(q)$.
In fact,  there exists for each dominant weight $\mu$
a unique (up to isomorphism) indecomposable tilting module $T_\mu$ with highest
weight  equal to $\mu$ (see [A, Cor.  2.6]). 
If $p$ is a primitive idempotent 
in the $k$-row quotient of $H_n(q)$ whose image in the maximal semisimple quotient
is a minimal idempotent ``belonging" to the simple module  $D^\mu$,
then $T_\mu \cong p V^{\otimes n}$. Here we have
identified the Young diagram $\mu$ with a dominant
weight of $sl_k$ in the usual way.

Tilting modules have filtrations by {\em standard} or {\em Weyl} modules
$\Delta_\la$, whose characters are given by the Weyl character formula. 
Consequently, the  character $\chi^\mu_T$ of the tilting module 
$T_\mu$ is  a linear combination of characters of Weyl modules, 
\begin{equation} \label{tilting characters}
\chi^\mu_T=\sum n_{\la\mu}\chi^\la ,
\end{equation}
where $n_{\la\mu}$ are non-negative integers.

It is well known that the multiplicities of Weyl modules in indecomposable
tilting modules coincides with the multiplicities of simple Hecke algebra
modules in Specht modules:
\begin{equation} \label{equal decomposition numbers}
n_{\la\mu} = d_{\la\mu}.
\end{equation}
This equality has been derived using a duality theory for quasi-heriditary algebras;
see [Do2], Chapter 4, and also [Do1], Lemma 3.1.  The equality can also be
derived quite directly from Schur-Weyl duality and an evaluation 
argument; we give this argument in section 6.2.

In [GW], we have extended the equality \ref{equal decomposition numbers}
to certain polynomial analogues of the decomposition numbers:
 \begin{equation} \label{equal decomposition polynomials}
n_{\la\mu}(v) = d_{\la\mu}(v).
\end{equation} 
Here the $n_{\la\mu}(v)$ are affine Kazhdan-Lusztig polynomials which satisfy
$n_{\la\mu}(1) = n_{\la\mu}$ ([S1], [S2]), and  the
$d_{\la\mu}(v)$ are polynomials associated to a
$U_v(\widehat{sl}_l)$-module, which satisfy 
$d_{\la\mu}(1) = d_{\la\mu}$ ([LLT], [Ar]).

\v\ni 
\subsection{}  We provide  a proof of the equality (\ref{equal
decomposition numbers}) by a fairly elementary evaluation argument.
Let $F$ be any field of characteristic 0, and  $q$ an element of $  F^\times$.
Consider the commuting actions of $H_n(q, F)$ and of the quantum group
$U_q(sl_k(F))$ on $(F^k)^{\otimes n}$.  (If $q$ is a root of unity, take the
Lusztig version of the quantum group.)  The standard basis of
 $(F^k)^{\otimes n}$ is a basis of weight vectors for the Cartan subalgebra
of $U_q(sl_k(F))$.   For each dominant integral weight $\gamma$, the
weight space  $W(\gamma) \subseteq  (F^k)^{\otimes n}$ is an $H_n(q, F)$
submodule.  For any idempotent $p \in H_n(q, F)$, $\Tr_{W(\gamma)}(p)$
is the dimension of $W(\gamma) \cap p (F^k)^{\otimes n}$, that is, the multiplicity
of the weight $\gamma$ in the $U_q(sl_k(F))$ module  $ p(F^k)^{\otimes n}$.

Now consider a field $K$ of characteristic 0, put $V(x) = K(x)^k$, and 
$V = K^k$.  Fix a primitive $l$-th root of unity $q \in K$.
Consider
the commuting actions of $H_n(x, K(x))$ and $U_x(sl_k)$ on $V(x)^{\otimes n}$.
For each dominant integral weight $\gamma$, the weight space
$W(\gamma, x)$ of $\gamma$ in $V(x)^{\otimes n}$ is an evaluable $H_n(x, K(x))$
module (with respect to the standard basis), and its ``restriction" to 
$H_n(q, K)$ is the weight space $W(\gamma)$ of $\gamma$ in $V^{\otimes n}$.
If $p \in H_n(x, K(x))$ is an evaluable idempotent, then according to 
Proposition \ref{Proposition 0.6.1}(a), 
$\Tr_{W(\gamma, x)}(p) = \Tr_{W(\gamma)}(p(q))$.  That is, the multiplicities  of 
the weight $\gamma$ in $p V(x)^{\otimes n}$ and in  $p(q) V^{\otimes n}$ are the
same.

Note also that if $e$ is a minimal idempotent in the semisimple algebra
$H_n(x, K(x))$ belonging to the simple module $S^\la$, then
$e V(x)^{\otimes n}$ is a simple $U_x(sl_k)$ with highest weight $\lambda$
and character $\chi^\la$ (by Schur-Weyl duality)  and 
$\Tr_{W(\gamma, x)}(e) = m_\gamma^\la$ is the multiplicity of $\gamma$ in
this module.

Finally, let $p \in H_n(x, K(x))$ be an evaluable idempotent such
that $p(q)$ is a primitive idempotent in $H_n(q, K)$ whose image in the
maximal semisimple quotient is minimal and belongs to the simple module $D^\mu$.
As observed above, $p(q) V^{\otimes n} \cong T_\mu$.
If $p$ is written as an orthogonal sum of minimal idempotents in $H_n(x, K(x))$, then
$d_{\la \mu}$ of these minimal idempotents satisfy $e S^\la \ne 0$.
Thus the multiplicity of  $\gamma$ in $T_\mu$ is
\begin{equation}
\Tr_{W(\gamma)}(p(q)) = \Tr_{W(\gamma, x)}(p) = \sum_\la d_{\la \mu} m_\gamma^\la.
\end{equation}
On the other hand, it follows from Equation (\ref{tilting characters}) that this
multiplicity is $\sum_\la n_{\la \mu} m_\gamma^\la$.  From the equations
$\sum_\la n_{\la \mu} m_\gamma^\la = \sum_\la d_{\la \mu} m_\gamma^\la$ for all
$\gamma$, one has $\sum_\la  n_{\la \mu} \chi^\la =  \sum_\la  d_{\la \mu} \chi^\la$,
and hence $ n_{\la \mu}  =  d_{\la \mu}$ for all $\la$, by linear independence of
the Schur functions.

\v\ni 
\subsection{}
We show  the connection between our Theorem
\ref{Theorem 4.2} and results obtained by Soergel: Let $\Pi$ be
the fundamental box, i.e the region of points $x$
in the  Weyl chamber for which $0<(\alpha_i,x)<l$ for any simple 
root $\alpha_i$, and  let $\Delta =\cup_{w\in S_n}\overline{w(\Pi)}$.
A slightly weaker version of our Theorem \ref{Theorem 4.2}
 can be rephrased as  follows:
\v\ni
\begin{theorem}\label{Theorem 6.2}
Let $\mu+\rho\in c+\Pi$, where $c$ is a  critical point
in the Weyl chamber such that also $c+\Delta$ is contained in the
closure of the Weyl chamber. Then
the multiplicities $n_{\la\mu}$ and $d_{\la\mu}$ are nonzero only
if $\la+\rho\in c+\Delta$.
\end{theorem}
\v
This
result  follows from Theorem 5.3 and Proposition 4.19 in [S1].
The latter also contains much more precise information about these 
multiplicities.

\v\ni 
\subsection{}
If $F$ is a field of characteristic $p$, one can find
an embedding of the $k$-row quotient of $FS_n$, for appropriate $n$,
into the $k$-row
quotient of a much larger symmetric group, using paths
in the region below $\la_1-\la_k=p^2$. One only needs to replace
the $q$-numbers by ordinary numbers (mod $p$) in all our proofs.
It should be possible to extend this embedding for arbitrary
$n$ also in the region above $\la_1-\la_k=p^2$,
although one certainly would have to modify parts of our
proof for the general case.

It should be remarked that a similar statement as Theorem \ref{Theorem 6.2}
does NOT hold for
tilting modules of algebraic groups over fields with positive
characteristic. In this case, there  are no uniform bounds for
the number of  irreducible generic characters which appear in the
character of an indecomposable tilting module.

\v\ni
\subsection{}
In the setting of algebraic groups, there is an interesting
connection between $Sl(k,F)$ Weyl modules $V_{\la}$ and $V_{p\la +(p-1)\rho}$
involving the Frobenius homomorphism (see [Ja, Section 10] for details).
One would expect a similar connection between $sl_k$-modules $V_\la$ and
$U_qsl_k$-modules  $V_{l\la +(l-1)\rho}$ for $q$ a primitive $l$-th
root of unity. More precisely,
it is easy to check 
that we have the following isomorphism of vector spaces:
$$T_{l\la +(l-1)\rho} \cong V_\la\otimes T_{(l-1)\rho},\eqno(*)$$
where $V_\la$ is an irreducible $sl_k$ module  with highest weight
$\la$. This follows easily from Weyl's dimension formula
(observe that the tilting modules in $(*)$ coincide with the
Weyl module with the same highest weight).  Comparison
with the setting for algebraic groups (see [Ja, 10.5]) 
would suggest that this can
be done in such a way that $V_\la$ is an $sl_k$-module, and
$T_{(l-1)\rho}$ is a $K_l$ module, where $K_l$ is the kernel of
Lusztig's Frobenius type homomorphism (see [Lu]). 
 It would seem
that our result Theorem \ref{Theorem 3.3.1} is in some sense dual to Lusztig's
result, with `dual' in the sense of Schur-Weyl duality.

Similarly, one would expect from the occurrence of  
a representation of the symmetric group within the Hecke algebra
a tensor action $\hat\otimes$ of  $Rep(sl_k)$ on direct sums of 
modules of the form $T_{l\la +(l-1)\rho}$ with certain categorical
properties. Using the vector space
isomorphism in $(*)$, one would expect a decomposition 
$$V_\mu\hat\otimes T_{l\la+(l-1)\rho}=\bigoplus_\nu c_{\la\mu}^\nu 
T_{l\nu +(l-1)\rho},$$
where $V_\mu , V_\nu$ are simple highest weight modules of $sl_k$
and $c_{\la\mu}^\nu$ is the multiplicity of $V_\nu$ in $V_\la\otimes V_\mu$.

\v\ni \subsection{}
It is  comparatively straightforward to prove estimates as the ones
in Theorem \ref{Theorem 6.2} for type $A$ also
for the multiplicities $n_{\la\mu}$ for  representations
of $U_qsp_{2k}$, or  for those representations of $U_qso_{2k+1}$ which appear
in some tensor power of the defining $k$-dimensional representation,
using a $q$-version of Brauer's centralizer algebras.

\v

\end{document}